\documentclass[11pt]{article}
\usepackage{latexsym}
\usepackage{amsmath,amsthm,amsfonts,amssymb,graphicx,epsfig,latexsym,color,mathrsfs}
\usepackage{epsf,enumerate}
\usepackage{fullpage}
\usepackage{tikz}
\usepackage{subcaption}
\usetikzlibrary{shapes,snakes}

\newtheorem{thm}{Theorem}
\newtheorem{lemma}[thm]{Lemma}

\newtheorem{prop}[thm]{Proposition}
\newtheorem{defn}[thm]{Definition}
{}

\theoremstyle{remark}

\newtheorem*{definition*}{Definition}
\newtheorem*{remark*}{Remark}

\def\qed{\hfill \ifhmode\unskip\nobreak\fi\quad\ifmmode\Box\else$\Box$\fi\\ }

\renewcommand{\phi}{\varphi}

\begin{document}

\title{Packing chromatic number of subdivisions of cubic graphs}
\date{\today}
\author{
J\' ozsef Balogh~\thanks{Department of Mathematics, University of Illinois at Urbana--Champaign, IL, USA and Moscow Institute of Physics and Technology, 9 Institutskiy per., Dolgoprodny, Moscow Region, 141701, Russian Federation, jobal@illinois.edu. 
Research of this author is partially supported by NSF Grant DMS-1500121 and by the Langan Scholar Fund
(UIUC).}
\and Alexandr Kostochka \thanks{Department of Mathematics, University of Illinois at Urbana--Champaign, IL, USA and
Sobolev Institute of Mathematics, Novosibirsk 630090, Russia, kostochk@math.uiuc.edu. Research of this author is supported in part by NSF grant
 DMS-1600592 and by grants 18-01-00353A  and 16-01-00499 of the Russian Foundation for Basic Research.}
 \and Xujun Liu\thanks{Department of Mathematics, University of Illinois at Urbana--Champaign, IL, USA,  xliu150@illinois.edu.}
 }

\maketitle

\begin{abstract}
A {\em packing $k$-coloring} of a graph $G$ is a partition of $V(G)$ into sets $V_1,\ldots,V_k$ such that for each $1\leq i\leq k$
the distance between any two distinct $x,y\in V_i$ is at least $i+1$. The {\em packing chromatic number}, $\chi_p(G)$, of a graph $G$ is the minimum $k$ such that $G$ has a packing $k$-coloring. For a graph $G$, let $D(G)$ denote the graph obtained from $G$ by subdividing every edge. 
The questions on the value of the maximum of $\chi_p(G)$  and of $\chi_p(D(G))$ over the class of subcubic graphs $G$ 
appear in several papers. Gastineau and  Togni asked whether $\chi_p(D(G))\leq 5$ for any subcubic $G$, and
later Bre\v sar,  Klav\v zar,  Rall and  Wash conjectured this, but no upper bound was proved.
 Recently the authors proved that $\chi_p(G)$  is not bounded in the class of subcubic graphs $G$. 
In contrast, in this paper we show that 
$\chi_p(D(G))$ is bounded in this class, and does not exceed $8$.

\end{abstract}

	\section{Introduction}	
For a positive integer $i$,  a set $S$ of vertices in a graph $G$ is {\em $\; i$-independent } if the distance in $G$ between any two
distinct vertices of $S$ is at least $i+1$. In particular, a $1$-independent set is simply an independent set.

A {\em packing $k$-coloring} of a graph $G$ is a partition of $V(G)$ into sets $V_1,\ldots,V_k$ such that for each $1\leq i\leq k$,
the set $V_i$ is $i$-independent. The {\em packing chromatic number}, $\chi_p(G)$, of a graph $G$, is the minimum $k$ such that $G$ has a packing $k$-coloring.
 The notion of packing $k$-coloring was introduced in 2008 by
Goddard,  Hedetniemi, Hedetniemi,  Harris and  Rall~\cite{GHHHR}
(under the name {\em broadcast coloring}) motivated by frequency assignment problems in broadcast networks.
The concept has attracted a considerable attention recently: there are around 30 papers on the topic 
(see e.g.~\cite{ANT1,BF,BKR1,BKR2,BKRW1,BKRW2,BKRW3,CJ1,FG1,FKL1,G1,GT1,S1} and references in them). In particular, 
Fiala and Golovach~\cite{FG1} proved
that finding the  packing chromatic number of a graph 
is NP-hard even in the class of trees. Sloper~\cite{S1} showed that there are graphs with maximum degree $4$ and arbitrarily large
packing chromatic number. In particular, coloring of {\em graph subdivisions} were considered.
 For a graph $G$, let $D(G)$ denote the graph obtained from $G$ by subdividing every edge. 

The questions on how large can $\chi_p(G)$ and $\chi_p(D(G))$ be if $G$ is a subcubic graph
(i.e., a graph with maximum degree at most $3$) were discussed in several papers (see~\cite{BKRW1,BKRW2,GT1,LBS1,S1}).
In particular, Gastineau and Togni~\cite{GT1} asked whether $\chi_p(D(G))\leq 5$ for every subcubic graph  $G$.
Bre\v sar, Klav\v zar, Rall, and Wash~\cite{BKRW2} later conjectured this and proved the validity of their conjecture for some
special classes of subcubic graphs (e.g., the class of generalized Petersen graph).
However, no upper bounds for the whole class of (sub)cubic graphs
were proved in either case. Recently, the authors~\cite{BKL} showed that $\chi_p(G)$ is not bounded in 
the class of cubic graphs and that `many' cubic graphs have
`high' packing chromatic number.

In contrast, in this paper we give the first upper bound on $\chi_p(D(G))$ for subcubic $G$: we
show that $\chi_p(D(G))$ is  bounded by $8$ in 
this class. We will prove the following slightly stronger result.

\begin{thm}\label{t1}
For every connected subcubic graph $G$, the graph $D(G)$ has a packing $8$-coloring such that color $8$ is used at most once.
\end{thm}

The theorem will be proved in the language of $S$-colorings introduced in~\cite{GX1} and used in~\cite{GT1,GX2}.

\begin{defn}
For a non-decreasing sequence $S=(s_1,s_2,\ldots,s_k)$ of positive integers, an $S${\em-coloring} of a graph $G$ is a partition of $V(G)$ into sets $V_1,\ldots,V_k$ such that for each $1 \le i \le k$ the distance between any two distinct $x,y \in V_i$ is at least $s_i+1$. 
\end{defn}

In particular, a $(1,\ldots,1)$-coloring is an ordinary coloring, and a $(1,2,\ldots,k)$-coloring is a packing $k$-coloring. For subcubic graphs, Gastineau and Togni~\cite{GT1} proved that they are $(1,1,2,2,2)$-colorable and $(1,2,2,2,2,2,2)$-colorable. We will use the
following observation of Gastineau and Togni~\cite{GT1}.

\begin{prop}[\cite{GT1} Proposition 1]\label{extension}
Let $G$ be a graph and $S
=(s_1, \ldots,s_k)$ be a non-decreasing sequence of integers. If $G$ is $S$-colorable then $D(G)$ is $(1,2s_1+1, \ldots,2s_k+1)$-colorable.
\end{prop}

In particular, if $G$ is $(1,1,2,2,3,3)$-colorable, then $D(G)$ has a packing $7$-coloring.
In view of this, by a
 {\em feasible} coloring of $G$ we call a
  coloring of $G$ with colors $1_a,1_b,2_a,2_b,3_a,3_b$ such that the distance between any two distinct vertices of color $i_x$
is at least $i+1$ for all $1\leq i\leq 3$ and $x\in \{a,b\}$.

\begin{defn}
A {\em k-degenerate graph} is a graph in which every subgraph has a vertex of degree at most $k$.
\end{defn}

In the next two sections  we discuss feasible coloring
 of $2$-degenerate subcubic graphs. In Section~2, we will 
 show that
 if a $2$-degenerate subcubic graph $G$ has a feasible coloring $f$ and $v,u$ are  vertices of $G$ with degree at most $2$, then we can
 change $f$ to another feasible coloring with some control on the colors of $v$ and $u$. The long proof of one of the lemmas, Lemma~\ref{c2},
 is postponed till the last section.
 Based on the lemmas of Section~2, in Section~3 we prove
 the following theorem (that gives a better bound than Theorem~\ref{t1} but for a more restricted class of graphs).
 
\begin{thm}\label{t2}
Every $2$-degenerate subcubic graph $G$ has a feasible coloring. In particular,
 $D(G)$ has a packing $7$-coloring.
\end{thm}
 
 In Section~4 we use Theorem~\ref{t2} and the lemmas in Section~2 to derive Theorem~\ref{t1}.
 In the final section we present a proof of Lemma~\ref{c2}.



\section{Lemmas on feasible coloring}

\begin{defn} For a positive integer $s$ and a vertex $a$ in a graph $G$, the {\em ball $B_{G}(a,s)$ in $G$ of radius $s$ with center $a$}
is $\{v\in V(G)\,:\; d_G(v,a)\leq s\}$, where $d_G(v,a)$ denotes the distance in $G$ between $v$ and $a$. We abbreviate $B_G(a,s)$ to $B(a,s)$ when the graph $G$ is clear from the context.
\end{defn}

\begin{defn}
For a positive integer $k$, a {\em $k$-vertex} is a vertex of degree exactly $k$.
\end{defn}

For $A=\{a_1, \ldots, a_n\} \subseteq V(G)$ and a coloring $f$, by $f(A)$ we mean $\{f(a_1), \ldots, f(a_n)\}$.

\begin{lemma}\label{c1}
Let $G$ be a subcubic graph and $f$ be a feasible coloring of $G$. Suppose there are $2$-vertices $u,v\in V(G)$ with
$f(u)=f(v)=2_a$. Let $N(u)=\{u_1,u_2\}$ and  $N(v)=\{v_1,v_2\}$. Then $G$ has a feasible coloring $g$ satisfying one of the following:\\
(a) $g(u)=2_a$ and $g(v)\in \{1_a,1_b\}$ or  $g(v)=2_a$ and $g(u)\in \{1_a,1_b\}$;\\
(b) $\{g(u),g(v)\}=\{2_a,2_b\}$;\\
(c) $\{g(u_1),g(u_2)\}=\{g(v_1),g(v_2)\}=\{1_a,1_b\}$, and exactly one of $u,v$ has color $2_a$.
\end{lemma}

{\bf Proof.} If $\{f(u_1),f(u_2)\}\neq \{1_a,1_b\}$, then we recolor $u$ with a color $\alpha\in \{1_a,1_b\}-\{f(u_1),f(u_2)\}$,
and $(a)$ holds. Thus by the symmetry between $u$ and $v$ we may assume 
\begin{equation}\label{1ab}
f(u_1)=f(v_1)=1_a\quad  \mbox{and}\quad  f(u_2)=f(v_2)=1_b.
\end{equation}
 Since $f(u)=f(v)=2_a$, $N(u)\cap N(v)=\emptyset$. In other words,
 \begin{equation}\label{distinct}
  \mbox{\em all vertices $u_1,u_2,v_1$ and $v_2$ are distinct.}
\end{equation}
Let $G_1$ denote the subgraph of $G$ induced by the vertices of colors $1_a$ and $1_b$. If $u_1$ and $u_2$ are in distinct
components of $G_1$, then after switching the colors in the component of $G_1$ containing $u_2$, we obtain a coloring
contradicting~\eqref{1ab}. Thus we may assume
\begin{equation}\label{1ab'}
  \mbox{\em $G$ has a $1_a,1_b$-colored $u_1,u_2$-path $P_u$ and a $1_a,1_b$-colored $v_1,v_2$-path $P_v$.}
\end{equation}

{\bf Case 1:} $u_1u_2\in E(G)$. If $|N(u_1)|=3$, then let $u_3\in N(u_1)-\{u,u_2\}$. Similarly, if
$|N(u_2)|=3$, then let $u_4\in N(u_2)-\{u,u_1\}$. If $2_b\notin f(N(u_1)\cup N(u_2))$, then after recoloring $u$ with $2_b$ we get a
coloring satisfying (b). Thus we may assume 
\begin{equation}\label{1ab''}
  \mbox{\em
$|N(u_1)|=3$ and $f(u_3)=2_b$.}
\end{equation}
Let $N(u_3)\subseteq \{u_1,u_5,u_6\}$.  If $2_a\notin f(N(u_3))$, then since $f(u_4)\neq 2_a$ (because
$d(u,u_4)=2$) after switching the colors of $u$ and $u_1$ 
we obtain a coloring satisfying $(a)$. So we may assume $f(u_5)=2_a$.

\begin{figure}[ht]\label{f1}
\hspace{5mm}
\begin{minipage}[b]{0.3\textwidth}
\begin{tikzpicture}[scale=0.5, transform shape]

\node[circle, draw=white!0, inner sep=0pt, minimum size=25pt, font=\huge] (11) at (1.15,0.60) {$2_a$};
\node[circle, draw=white!0, inner sep=0pt, minimum size=25pt, font=\huge] (12) at (-1.88,-1.5) {$1_a$};
\node[circle, draw=white!0, inner sep=0pt, minimum size=25pt, font=\huge] (13) at (7.3,0.75) {$2_a$};
\node[circle, draw=white!0, inner sep=0pt, minimum size=25pt, font=\huge] (14) at (2.7,-1.4) {$1_b$};
\node[circle, draw=white!0, inner sep=0pt, minimum size=25pt, font=\huge] (15) at (4.2,-1.2) {$1_a$};
\node[circle, draw=white!0, inner sep=0pt, minimum size=25pt, font=\huge] (16) at (8.7,-1.2) {$1_b$};
\node[circle, draw=white!0, inner sep=0pt, minimum size=25pt, font=\huge] (17) at (-1.88,-3.5) {$2_b$};
\node[circle, draw=white!0, inner sep=0pt, minimum size=25pt, font=\huge] (18) at (0.5,-5.3) {$1_b$};
\node[circle, draw=white!0, inner sep=0pt, minimum size=25pt, font=\huge] (19) at (-2.8,-5.3) {$2_a$};

\node[circle, draw=black!80,  inner sep=0pt, minimum size=25pt, font=\huge] (1) at (0.5,0) {$u$};
\node[circle, draw=black!80,  inner sep=0pt, minimum size=25pt, font=\huge] (3) at (-1,-2) { $u_1$};
\node[circle, draw=black!80,  inner sep=0pt, minimum size=25pt, font=\huge] (4) at (2,-2) {$u_2$};
\node[circle, draw=black!80,  inner sep=0pt, minimum size=25pt, font=\huge] (5) at (-1,-4) {$u_3$};
\node[circle, draw=black!80,  inner sep=0pt, minimum size=25pt, font=\huge] (6) at (-2,-6) {$u_5$};
\node[circle, draw=black!80,  inner sep=0pt, minimum size=25pt, font=\huge] (7) at (0,-6) {$u_6$};
\node[circle, draw=black!80,  inner sep=0pt, minimum size=25pt, font=\huge] (8) at (2,-4) {$u_4$};
\node[circle, draw=black!80,  inner sep=0pt, minimum size=25pt, font=\huge] (2) at (6.5,0) {$v$};
\node[circle, draw=black!80,  inner sep=0pt, minimum size=25pt, font=\huge] (9) at (5,-2) {$v_1$};
\node[circle, draw=black!80,  inner sep=0pt, minimum size=25pt, font=\huge] (10) at (8,-2) {$v_2$};
\draw  (1) edge (3);
\draw  (1) edge (4);
\draw  (3) edge (4);
\draw  (3) edge (5);
\draw  (4) edge (8);
\draw  (5) edge (6);
\draw  (5) edge (7);
\draw  (2) edge (9);
\draw  (2) edge (10);

\end{tikzpicture}

\caption{Case 1.1.}
\label{case 1.1.}
\end{minipage}
\hspace{15mm}
\begin{minipage}[b]{0.4\textwidth}
\begin{tikzpicture}[scale=0.4, transform shape]
\node[circle, draw=white!0,  inner sep=0pt, minimum size=25pt, font=\huge] (11) at (0,1) {$2_a$};
\node[circle, draw=white!0,  inner sep=0pt, minimum size=25pt, font=\huge] (12) at (-4,-1) {$1_a$};
\node[circle, draw=white!0,  inner sep=0pt, minimum size=25pt, font=\huge] (14) at (4,-1) {$1_b$};

\node[circle, draw=white!0,  inner sep=0pt, minimum size=25pt, font=\huge] (15) at (11,1) {$2_a$};

\node[circle, draw=white!0,  inner sep=0pt, minimum size=25pt, font=\huge] (19) at (-7,-7.3) {$1_a$};

\node[circle, draw=white!0,  inner sep=0pt, minimum size=25pt, font=\huge] (20) at (-5,-7.3) {$2_a$};
\node[circle, draw=white!0,  inner sep=0pt, minimum size=25pt, font=\huge] (40) at (-5,-8.2) {$2_b$};
\node[circle, draw=white!0, inner sep=0pt, minimum size=25pt, font=\huge] (21) at (8.5,-7.5) {$A$};
\node[circle, draw=white!0, inner sep=0pt, minimum size=25pt, font=\huge] (22) at (-2,-3) {$2_b$};

\node[circle, draw=white!0, inner sep=0pt, minimum size=25pt, font=\huge] (23) at (-3,-7.3) {$1_b$};
\node[circle, draw=white!0, inner sep=0pt, minimum size=25pt, font=\huge] (24) at (2,-3) {$1_a$};

\node[circle, draw=white!0, inner sep=0pt, minimum size=25pt, font=\huge] (26) at (13.5,-1) {$1_b$};

\node[circle, draw=white!0, inner sep=0pt, minimum size=25pt, font=\huge] (30) at (8.5,-1) {$1_a$};

\node[circle, draw=white!0, inner sep=0pt, minimum size=25pt, font=\huge] (32) at (-1,-7.3) {$1_a$};
\node[circle, draw=white!0, inner sep=0pt, minimum size=25pt, font=\huge] (41) at (-1,-8.2) {$2_a$};

\node[circle, draw=white!0, inner sep=0pt, minimum size=25pt, font=\huge] (37) at (-6,-3) {$1_b$};

\node[circle, draw=white!0, inner sep=0pt, minimum size=25pt, font=\huge] (39) at (1,-7.3) {$1_b$};

\node[circle, draw=black!80, inner sep=0pt, minimum size=25pt, font=\huge] (1) at (0,0) {$u$};
\node[circle, draw=black!80, inner sep=0pt, minimum size=25pt, font=\huge] (3) at (-4,-2) { $u_1$};
\node[circle, draw=black!80, inner sep=0pt, minimum size=25pt, font=\huge] (4) at (4,-2) {$u_2$};
\node[circle, draw=black!80, inner sep=0pt, minimum size=25pt, font=\huge] (5) at (-6,-4) {$u_3$};
\node[circle, draw=black!80, inner sep=0pt, minimum size=36pt, font=\huge] (6) at (-7,-6) {$u_7$};
\node[circle, draw=black!80, inner sep=0pt, minimum size=36pt, font=\huge] (7) at (-5,-6) {$u_9$};
\node[circle, draw=black!80, inner sep=0pt, minimum size=25pt, font=\huge] (8) at (2,-4) {$u_4$};
\node[circle, draw=black!80, inner sep=0pt, minimum size=25pt, font=\huge] (2) at (11,0) {$v$};
\node[circle, draw=black!80, inner sep=0pt, minimum size=25pt, font=\huge] (9) at (8.5,-2) {$v_1$};
\node[circle, draw=black!80, inner sep=0pt, minimum size=25pt, font=\huge] (10) at (13.5,-2) {$v_2$};

\node[circle, draw=black!80, inner sep=0pt, minimum size=36pt, font=\huge] (11) at (1,-6) {$u_8$};
\node[circle, draw=black!80, inner sep=0pt, minimum size=36pt, font=\huge] (12) at (7,-6) {$u_{14}$};
\node[circle, draw=black!80, inner sep=0pt, minimum size=36pt, font=\huge] (13) at (5,-6) {$u_{12}$};
\node[circle, draw=black!80, inner sep=0pt, minimum size=36pt, font=\huge] (14) at (3,-6) {$u_{10}$};

\node[circle, draw=black!80, inner sep=0pt, minimum size=25pt, font=\huge] (29) at (-2,-4) {$u_5$};
\node[circle, draw=black!80, inner sep=0pt, minimum size=25pt, font=\huge] (v7) at (6,-4) {$u_6$};

\node[circle, draw=black!80, inner sep=0pt, minimum size=36pt, font=\huge] (v1) at (-3,-6) {$u_{11}$};
\node[circle, draw=black!80, inner sep=0pt, minimum size=30pt, font=\huge] (v2) at (-1,-6) {$u_{13}$};

\draw  (2) edge (9);
\draw  (2) edge (10);
\draw  (1) edge (3);
\draw  (1) edge (4);
\draw  (4) edge (8);
\draw  (4) edge (v7);
\draw  (29) edge (3);
\draw  (5) edge (3);
\draw  (5) edge (6);
\draw  (5) edge (7);
\draw  (v1) edge (29);
\draw  (29) edge (v2);
\draw  (8) edge (11);
\draw  (8) edge (14);
\draw  (v7) edge (13);
\draw  (v7) edge (12);
\draw  (5.6,-5.5) ellipse (3.5 and 2.5);
\end{tikzpicture}
\caption{Case 2.1.}
\label{case 2.1.}
\end{minipage}
\end{figure}


{\bf Case 1.1:} $|N(u_2)|<3$ or $f(u_4)\neq 2_b$.
 If $1_b\notin f(N(u_3))$, then we can recolor $u_3$ with $1_b$. By the case, we can recolor $u$ with $2_b$ to obtain a coloring satisfying $(b)$. So we may assume $f(u_6)=1_b$ (See Figure~\ref{case 1.1.}).
 Then the coloring $g$
obtained from $f$ by recoloring  $u$ and $u_3$ with $1_a$ and $u_1$ with $2_b$ satisfies $(a)$.

{\bf Case 1.2:} $|N(u_2)|=3$ and $f(u_4)= 2_b$. If $u_4=u_3$, then $N(u_3)=\{u_1,u_2,u_5\}$. Then $u$ has no vertices of
color $3_a$ at distance at most $3$, so after recoloring $u$ with $3_a$, we obtain a coloring $g$ satisfying (c).
Thus, $u_4\neq u_3$. 

{\bf Case 1.2.1:} $1_b\notin f(N(u_3))$. We recolor $u_3$ with $1_b$. If $2_a \notin f(N(u_4)-u_2)$, then we recolor $u_2$ with $2_a$ and $u$ with $1_b$ to obtain a coloring satisfying $(a)$. If $1_a \notin f(N(u_4)-u_2)$, then we recolor $u_4$ with $1_a$, $u_2$ with $2_b$, and $u$ with $1_b$ to obtain a coloring satisfying $(a)$. Thus, we may assume 

$$
 f(N(u_4)-u_2)=\{1_a,2_a\}.
$$

Then recoloring $u_4$ with $1_b$, $u_2$ with $2_b$, and $u$ with $1_b$, we obtain a coloring satisfying $(a)$.

{\bf Case 1.2.2:} $1_b\in f(N(u_3))$. Since $f(u_5)=2_a$, this means $u_6$ exists and
 $f(u_6)=1_b$. Then we recolor $u_3$ and $u_2$ with $1_a$ and $u_1$ with $1_b$. If $2_a \notin f(N(u_4)-u_2)$, then we recolor $u_2$ with $2_a$ and $u$ with $1_a$ to obtain a coloring satisfying $(a)$. If $1_b \notin f(N(u_4)-u_2)$, then we recolor $u_4$ with $1_b$ and $u$ with $2_b$ to obtain a coloring satisfying $(b)$. Thus, we may assume 

$$
 f(N(u_4)-u_2)=\{1_b,2_a\}.
$$
 
Then we recolor $u_4$ with $1_a$, $u_2$ with $2_b$, and $u$ with $1_a$ to obtain a coloring satisfying $(a)$.

{\bf Case 2:} $u_1u_2\notin E(G)$. Then we may assume that $N(u_1)\subseteq \{u,u_3,u_5\}$,
$N(u_2)\subseteq \{u,u_4,u_6\}$ and
 by~\eqref{1ab'},  $f(u_3)=1_b$ and $f(u_4)=1_a$. Furthermore, since by the case, $u_3\neq u_2$, we may assume
 that $N(u_3)\subseteq \{u_1,u_7,u_9\}$ and $f(u_7)=1_a$. It is possible that $u_7=u_4$, but this will not affect the proof below.
 Similarly, we will assume
 that $N(u_4)\subseteq \{u_2,u_8,u_{10}\}$ and $f(u_8)=1_b$. As in Case 1, $2_b\in f(N(u_1)\cup N(u_2)),$ since otherwise we can recolor $u$ with $2_b$ and (b) will hold. In our notation,
 this means $2_b\in \{f(u_5),f(u_6)\}$. By symmetry, we will assume $f(u_5)=2_b$. We also will assume
 $N(u_5)\subseteq \{u_1,u_{11},u_{13}\}$ and $N(u_6)\subseteq \{u_2,u_{12},u_{14}\}$, where some vertices can coincide.

{\bf Case 2.1:} $|N(u_2)|<3$ or $f(u_{6})\neq 2_b$.  If $1_b\notin f(N(u_5))$, then we can recolor $u_5$ with $1_b$, 
and then $u$ with $2_b$. The resulting coloring satisfies (b). So we may assume
  $f(u_{11})=1_b$. If $2_a\notin \{f(u_9),f(u_{13})\}$, then by switching
the colors of $u$ and $u_1$, we obtain a coloring satisfying $(a)$. Thus $2_a\in \{f(u_9),f(u_{13})\}$.
If $f(u_9)=2_a$ and $f(u_{13})\neq 1_a$ or if $f(u_{13})=2_a$ and $f(u_{9})\neq 2_b$, then after switching the colors of $u_1$ and $u_5$ and recoloring $u$ with $1_a$,
we again get a coloring satisfying $(a)$. So,
\begin{equation}\label{u9}
  \mbox{\em either 
  $f(u_9)=2_a$ and $f(u_{13})= 1_a$ or  $f(u_{13})=2_a$ and $f(u_{9})= 2_b$.}
\end{equation}
If $u_6$ does not exist, then by~\eqref{u9}, the only vertex in $B(u,3)-(N(u) \cup \{u\})$ that can be colored with $3_a$ or $3_b$
is $u_{10}$. Thus after recoloring $u$ with a color in $\{3_a,3_b\}-f(u_{10})$ we obtain a coloring satisfying (c). So suppose
$u_6$ exists.
Let  
$A=\{u_6,u_{10},u_{12},u_{14}\}$. If $1_a\notin \{f(u_{12}),f(u_{14})\}$, then we can recolor $u_6$ with $1_a$ without changing color of any other vertex. Thus we may assume
\begin{equation}\label{u6}
 1_a\in f(A).
\end{equation}
If a color $x\in \{2_a,2_b\}$ is not in $f(A)$, then after recoloring $u_2$ with $x$ and $u$ with $1_b$, we get a coloring satisfying $(a)$.
Thus
\begin{equation}\label{u6'}
 2_a,2_b\in f(A).
\end{equation}
By the argument above, in particular, by~\eqref{u9}, colors $3_a$ and $3_b$ are not used on vertices in\\
$B=\{u_1,u_2,u_3,u_4,u_5,u_7,u_8,u_9,u_{11},u_{13}\}$. If at least one of them, say $3_a$, is also not used on $A$, then
after recoloring $u$ with $3_a$, we obtain a coloring satisfying (c). Thus
\begin{equation}\label{u6''}
 3_a,3_b\in f(A) \text{ (See Figure~\ref{case 2.1.})}.
\end{equation}
Since $|f(A)| \le 4$, relations~\eqref{u6},~\eqref{u6'} and~\eqref{u6''} cannot hold at the same time, a contradiction.

{\bf Case 2.2:} $|N(u_2)|=3$ and $f(u_{6})= 2_b$. Suppose first that $u_{6}=u_{5}$ and that $N(u_5)=\{u_1,u_2,u_{11}\}$.
If  $f(u_9)\neq 2_b$ and $f(u_{11})\neq 1_a$, then  after switching the colors of $u_1$ and $u_5$ and recoloring $u$ with $1_a$,
we  get a coloring satisfying $(a)$. So,
  $f(u_9)=2_b$ or $f(u_{11})= 1_a$. Similarly, considering switching colors of $u_2$ and $u_5$, we obtain that $f(u_{10})=2_b$ or $f(u_{11})= 1_b$. Together, this means
\begin{equation}\label{u5}
  \mbox{\em the colors of at least two vertices in $\{u_9,u_{10},u_{11}\}$ are in $\{1_a,1_b,2_b\}$.}
\end{equation}
By~\eqref{u5}, some color $y\in \{3_a,3_b\}$ is not used on $B(u,3)$. Then after recoloring 
$u$ with $y$, we obtain a coloring satisfying (c).

Now we assume  $u_{6}\neq u_{5}$. If $1_a\notin \{f(u_{12}), f(u_{14})\}$, then after recoloring $u_6$ with $1_a$, we
get Case 2.1. Thus below we assume $f(u_{12})=1_a$. If $2_a \notin \{f(u_{10}), f(u_{14})\}$, then we obtain a coloring satisfying $(a)$ by switching the colors of $u$ and $u_2$. Thus, $2_a\in \{f(u_{10}),f(u_{14})\}$. If $f(u_{14}) \neq 1_b$ and $f(u_{10}) \neq 2_b$, then after switching the colors of $u_2$ and $u_6$ and recoloring $u$ with $1_b$,
we again get a coloring satisfying $(a)$. So,
\begin{equation}\label{u10}
  \mbox{\em either 
  $f(u_{10})=2_a$ and $f(u_{14})= 1_b$ or  $f(u_{10})=2_b$ and $f(u_{14})= 2_a$.}
\end{equation}
Let $A=\{u_9,u_{11},u_{13}\}$. If $2_a \notin f(A)$, then we obtain a coloring satisfying $(a)$ by switching the colors of $u$ and $u_1$. Thus, 
\begin{equation}\label{2afa}
2_a \in f(A).
\end{equation}
If $1_a \notin f(\{u_{11},u_{13}\})$ and $f(u_{9}) \neq 2_b$, then after switching the colors of $u_1$ and $u_5$ and recoloring $u$ with $1_a$,
we again get a coloring satisfying $(a)$. Therefore,
\begin{equation}\label{1a2bfa}
1_a \in f(\{u_{11},u_{13}\})\text{ or } f(u_{9}) = 2_b.
\end{equation}
 By the argument above, in particular, by~\eqref{u10}, colors $3_a$ and $3_b$ are not used on vertices in\\
$B=\{u_1,u_2,u_3,u_4,u_5,u_7,u_8,u_{10},u_{12},u_{14}\}$. If at least one of them, say $3_a$, is also not used on $A$, then
after recoloring $u$ with $3_a$, we obtain a coloring satisfying (c). Thus, 
\begin{equation}\label{3a3bfa}
 3_a,3_b\in f(A).
\end{equation}
Since $|f(A)| \le 3$, relations~\eqref{2afa},~\eqref{1a2bfa}, and~\eqref{3a3bfa} cannot hold at the same time, a contradiction. \hfill\qed

Our second lemma is:

\begin{lemma}\label{c2}
Let $G$ be a subcubic graph and $f$ be a feasible coloring of $G$. Suppose there is a $2$-vertex $u \in V(G)$ with $N(u)=\{u_1,u_2\}$. If $f(u) \in \{3_a,3_b\}$, then we can recolor some vertices of $G$ so that the resulting coloring $g$ is feasible and satisfies the following:\\
(a) $g(u) \notin \{3_a,3_b\}$, and\\
(b) at most one vertex is recolored into $3_a$ or $3_b$, and this vertex (if there is such a vertex) is
at distance at most $3$ from $u$ and has degree $3$ in $G$, and at most one vertex of $f$-color $3_a$ or $3_b$ apart from $u$ is recolored into some other color, and this vertex (if there is such a vertex) has new color in $\{1_a, 1_b\}$. 
\end{lemma}

The proof of this lemma is a long case analysis, so we postpone it to the last section.


\section{Proof of Theorem~\ref{t2}}

We prove the theorem by induction on the number $n$ of vertices. When $n \le 6,$ the claim holds obviously, since we have 6 colors. When $n>6$, we assume the argument holds for every graph with fewer than $n$ vertices. Let $G$ be any $2$-degenerate subcubic graph with $n$
vertices. We may assume $G$ is connected.
Since $G$ is $2$-degenerate, it has a vertex, say $w$, with degree at most 2. 

{\bf Case 1:}  $d(w)=1$. Let $N(w)=w'$. Since $G-w$ is an $(n-1)$-vertex connected subcubic graph with $d_{G-w}(w') \le 2$, by the induction hypothesis, $G-w$ has a $(1,1,2,2,3,3)$-coloring $f$. We color $w$ with a color $x \in \{1_a,1_b\}-f(w')$ to extend $f$ to $G$. 

{\bf Case 2:}
$d(w)=2$. Let $N(w)=\{w_1,w_2\}$. Note that $G-w$ has at most two connected components and each connected component is a connected $2$-degenerate subcubic graph with less than $n$ vertices. By the induction hypothesis, $G-w$ has a feasible coloring $f$. We may assume that $|N_{G-w}(w_1)|=|N_{G-w}(w_2)|=2$. Otherwise we can first apply the induction hypothesis to obtain a $(1,1,2,2,3,3)$-coloring $f$ on $G-w$, then add leaves (vertices of degree one) to $w_1$ and $w_2$ to obtain a new graph $G'$ with $|N_{G'-w}(w_1)|=|N_{G'-w}(w_2)|=2$, then assign proper colors to those leaves we just added to obtain a $(1,1,2,2,3,3)$-coloring $f'$ on $G'-w$, then prove that $G'$ has a $(1,1,2,2,3,3)$-coloring, which can be used to get our desired coloring on $G$. So below we assume
$N(w_1)=\{w,w_3,w_4\}$ and $N(w_2)=\{w,w_5,w_6\}$.

By Lemma~\ref{c2}, $G-w$ has a feasible coloring $f_1$ such that $f_1(w_1)\notin \{3_a,3_b\}$. Then by Lemma~\ref{c2}
again, $G-w$ also has a feasible coloring $f_2$ such that $f_2(w_2)\notin \{3_a,3_b\}$ and no vertex of degree $2$
in $G-w$ changed its color to $3_a$ or $3_b$. Thus we also have $f_2(w_1)\notin \{3_a,3_b\}$. Therefore, $G-w$ has a feasible coloring $f_2$ such that $f_2(w_1) \notin \{3_a,3_b\}$ and $f_2(w_2) \notin \{3_a,3_b\}$.

{\bf Case 2.1:} Either  $f_2(w_1)\neq f_2(w_2)$ or $f_2(w_1)= f_2(w_2)\in\{1_a,1_b\}$.
 If $\{f_2(w_1),f_2(w_2)\} \neq \{1_a,1_b\}$, then we extend $f_2$ to $G$ by assigning $f_2(w)=\alpha \in \{1_a,1_b\}-\{f_2(w_1),f_2(w_2)\}$. By the case, if $f_2(w_1)=f_2(w_2)$, then $f_2(w_1)=f_2(w_2) \in \{1_a,1_b\}$. Therefore, the extension of $f_2$ to $G$ is feasible since we do not introduce new conflicts between $w_1$ and $w_2$ by adding $w$. Thus, we may assume 
\begin{equation}\label{w1w2}
f_2(w_1)=1_a \quad  \mbox{and}\quad  f_2(w_2)=1_b.
\end{equation}
 If $w_1$ and $w_2$ are in distinct
components of the subgraph $G_2$ of $G-w$ induced by the vertices of colors $1_a$ and $1_b$ in $f_2$, 
then after switching the colors $1_a$ and $1_b$ with each other in the component of $G_2$ containing $w_2$, we obtain a coloring
contradicting~\eqref{w1w2}. Thus we may assume
\begin{equation}\label{1a1b'''}
  \mbox{\em $G-w$ has a $1_a,1_b$-colored $w_1,w_2$-path $P_w$.}
\end{equation}
In particular, we may assume $f_2(w_3)=1_b$ and $f_2(w_5)=1_a$
(possibly, $w_3=w_2$ and then $w_5=w_1$).

If $\{2_a,2_b\} \nsubseteq f_2(N(w_1) \cup N(w_2)-\{w\})$, then we can extend $f_2$ to $G$ by assigning $f_2(w)=\beta \in \{2_a,2_b\}-f_2(N(w_1) \cup N(w_2)-\{w\})$. Thus, we may assume  
\begin{equation}\label{w1w2'}
|N(w_1)|=|N(w_2)|=3, \{2_a,2_b\} \subseteq f_2(N(w_1) \cup N(w_2)-\{w\}), \text{ and by symmetry }
\end{equation}
\begin{equation}\label{w3w4}
f_2(w_4)=2_a \quad  \mbox{and}\quad   f_2(w_6)=2_b.
\end{equation}
If $1_b \notin f_2(N(w_4)-w_2)$, then we can extend $f_2$ to a feasible coloring of $G$ by recoloring $w_4$ with
$1_b$ and letting $f_2(w)=2_a$. By this and the symmetric statement for $w_6$ we can assume that
\begin{equation}\label{w3w4n}
\mbox{ \em $w_4$ has a neighbor $w_7$ with $f_2(w_7)=1_b$ and $w_6$ has a neighbor $w_8$ with $f_2(w_8)=1_a$.} 
\end{equation}

{\bf Case 2.1.1:} $w_1w_2 \in E(G)$ (i.e.,
$w_3=w_2$ and  $w_5=w_1$). 
If $1_a \notin f_2(N(w_4)-w_1)$, then we obtain a feasible coloring on $G$ by switching colors of $w_1$ and $w_4$, assigning $1_a$ to $w$, and using $f_2$ on other vertices. Therefore, by \eqref{w3w4n}, we may assume $f_2(N(w_4)-w_1)=\{1_a,1_b\}$. Similarly, by \eqref{w3w4n}, we may assume $f_2(N(w_6)-w_2)=\{1_a,1_b\}$ (See Figure~\ref{case 2.1.1.}). With \eqref{w1w2}, \eqref{w3w4}, and the case, $3_a \notin f_2(B(w,3)-\{w\})$ and we can extend $f_2$ to $G$ by assigning $f_2(w)=3_a.$



{\bf Case 2.1.2:} $w_1w_2 \notin E(G).$ If $N(w_3)\cup N(w_4)$ does not contain a vertex $w_9$ of color $2_b$, then we can
recolor $w_1$ with $2_b$ and color $w$ with $1_a$. So we may assume that $N(w_3)\cup N(w_4)$  contains a vertex $w_9$ of color $2_b$ and symmetrically $N(w_5)\cup N(w_6)$  contains a vertex $w_{10}$ of color $2_a$. Furthermore, if
$1_a \notin f_2(N(w_4)-w_1)$ and $2_a \notin f_2(N(w_3)-w_1)$, then we can recolor $w_1$ with $2_a$
and color $w$ and $w_4$ with $1_a$. With \eqref{1a1b'''} and \eqref{w3w4n}, all vertices in $B(w_1,2)-w$ have colors in $\{1_a,1_b,2_a,2_b\}$.
Symmetrically, we can assume all vertices in $B(w_2,2)-w$ have colors in $\{1_a,1_b,2_a,2_b\}$ (See Figure~\ref{case 2.1.2.}). Then we can color $w$ with $3_a$.

\begin{figure}[ht]\label{f8}
\hspace{15mm}
\begin{minipage}[b]{0.25\textwidth}
\begin{tikzpicture}[scale=0.6, transform shape]

\node[circle, draw=white!0, inner sep=0pt, minimum size=25pt, font=\huge] (11) at (0.5,1) {$3_a$};
\node[circle, draw=white!0, inner sep=0pt, minimum size=25pt, font=\huge] (12) at (-1.9,-1.5) {$1_a$};

\node[circle, draw=white!0, inner sep=0pt, minimum size=25pt, font=\huge] (14) at (2.8,-1.4) {$1_b$};

\node[circle, draw=white!0, inner sep=0pt, minimum size=25pt, font=\huge] (16) at (2.5,-6) {$1_b$};
\node[circle, draw=white!0, inner sep=0pt, minimum size=25pt, font=\huge] (17) at (2.8,-3.5) {$2_b$};
\node[circle, draw=white!0, inner sep=0pt, minimum size=25pt, font=\huge] (18) at (-0.5,-6) {$1_b$};
\node[circle, draw=white!0, inner sep=0pt, minimum size=25pt, font=\huge] (19) at (-1.9,-3.5) {$2_a$};
\node[circle, draw=white!0, inner sep=0pt, minimum size=25pt, font=\huge] (20) at (1.5,-6) {$1_a$};

\node[circle, draw=white!0, inner sep=0pt, minimum size=25pt, font=\huge] (22) at (-1.5,-6) {$1_a$};

\node[circle, draw=black!80,  inner sep=0pt, minimum size=25pt, font=\huge] (1) at (0.5,0) {$w$};
\node[circle, draw=black!80,  inner sep=0pt, minimum size=25pt, font=\huge] (3) at (-1,-2) { $w_1$};
\node[circle, draw=black!80,  inner sep=0pt, minimum size=25pt, font=\huge] (4) at (2,-2) {$w_2$};
\node[circle, draw=black!80,  inner sep=0pt, minimum size=25pt, font=\huge] (5) at (-1,-4) {$w_4$};

\node[circle, draw=black!80,  inner sep=0pt, minimum size=25pt, font=\huge] (8) at (2,-4) {$w_6$};

\draw  (1) edge (3);
\draw  (1) edge (4);
\draw  (3) edge (4);
\draw  (3) edge (5);
\draw  (4) edge (8);

\node (v1) at (-1.5,-5.5) {};
\node (v2) at (-0.5,-5.5) {};
\node (v3) at (1.5,-5.5) {};
\node (v4) at (2.5,-5.5) {};
\draw  (5) edge (v1);
\draw  (5) edge (v2);
\draw  (v3) edge (8);
\draw  (8) edge (v4);
\end{tikzpicture}

\caption{Case 2.1.1.}
\label{case 2.1.1.}
\end{minipage}
\hspace{15mm}
\begin{minipage}[b]{0.4\textwidth}
\begin{tikzpicture}[scale=0.5, transform shape]
\node[circle, draw=white!0,  inner sep=0pt, minimum size=25pt, font=\huge] (11) at (0,1) {$3_a$};
\node[circle, draw=white!0,  inner sep=0pt, minimum size=25pt, font=\huge] (12) at (-4,-1) {$1_a$};
\node[circle, draw=white!0,  inner sep=0pt, minimum size=25pt, font=\huge] (14) at (4,-1) {$1_b$};

\node[circle, draw=white!0,  inner sep=0pt, minimum size=25pt, font=\huge] (15) at (6,-3) {$2_b$};

\node[circle, draw=white!0,  inner sep=0pt, minimum size=25pt, font=\huge] (19) at (-6.5,-6.5) {$1_a$};

\node[circle, draw=white!0,  inner sep=0pt, minimum size=25pt, font=\huge] (20) at (-5.5,-6.5) {$2_a$};
\node[circle, draw=white!0,  inner sep=0pt, minimum size=25pt, font=\huge] (40) at (-2.5,-6.5) {$2_b$};

\node[circle, draw=white!0, inner sep=0pt, minimum size=25pt, font=\huge] (22) at (-2,-3) {$2_a$};

\node[circle, draw=white!0, inner sep=0pt, minimum size=25pt, font=\huge] (23) at (-5.5,-7.5) {$2_b$};
\node[circle, draw=white!0, inner sep=0pt, minimum size=25pt, font=\huge] (24) at (2,-3) {$1_a$};

\node[circle, draw=white!0, inner sep=0pt, minimum size=25pt, font=\huge] (26) at (1.5,-6.5) {$1_b$};

\node[circle, draw=white!0, inner sep=0pt, minimum size=25pt, font=\huge] (30) at (5.5,-6.5) {$1_a$};
\node[circle, draw=white!0, inner sep=0pt, minimum size=25pt, font=\huge] (42) at (2.5,-7.5) {$2_b$};
\node[circle, draw=white!0, inner sep=0pt, minimum size=25pt, font=\huge] (43) at (6.5,-7.5) {$2_a$};
\node[circle, draw=white!0, inner sep=0pt, minimum size=25pt, font=\huge] (44) at (6.5,-6.5) {$1_b$};

\node[circle, draw=white!0, inner sep=0pt, minimum size=25pt, font=\huge] (32) at (-2.5,-7.5) {$1_a$};
\node[circle, draw=white!0, inner sep=0pt, minimum size=25pt, font=\huge] (41) at (2.5,-6.5) {$2_a$};

\node[circle, draw=white!0, inner sep=0pt, minimum size=25pt, font=\huge] (37) at (-6,-3) {$1_b$};

\node[circle, draw=white!0, inner sep=0pt, minimum size=25pt, font=\huge] (39) at (-1.5,-6.5) {$1_b$};

\node[circle, draw=black!80, inner sep=0pt, minimum size=25pt, font=\huge] (1) at (0,0) {$w$};
\node[circle, draw=black!80, inner sep=0pt, minimum size=25pt, font=\huge] (3) at (-4,-2) { $w_1$};
\node[circle, draw=black!80, inner sep=0pt, minimum size=25pt, font=\huge] (4) at (4,-2) {$w_2$};
\node[circle, draw=black!80, inner sep=0pt, minimum size=25pt, font=\huge] (5) at (-6,-4) {$w_3$};

\node[circle, draw=black!80, inner sep=0pt, minimum size=25pt, font=\huge] (8) at (2,-4) {$w_5$};

\node[circle, draw=black!80, inner sep=0pt, minimum size=25pt, font=\huge] (29) at (-2,-4) {$w_4$};
\node[circle, draw=black!80, inner sep=0pt, minimum size=25pt, font=\huge] (v7) at (6,-4) {$w_6$};

\draw  (1) edge (3);
\draw  (1) edge (4);
\draw  (4) edge (8);
\draw  (4) edge (v7);
\draw  (29) edge (3);
\draw  (5) edge (3);

\node (v1) at (-6.5,-6) {};
\node (v2) at (-5.5,-6) {};
\node (v3) at (-2.5,-6) {};
\node (v4) at (-1.5,-6) {};
\node (v5) at (1.5,-6) {};
\node (v6) at (2.5,-6) {};
\node (v8) at (5.5,-6) {};
\node (v9) at (6.5,-6) {};
\draw  (5) edge (v1);
\draw  (5) edge (v2);
\draw  (29) edge (v3);
\draw  (29) edge (v4);
\draw  (8) edge (v5);
\draw  (8) edge (v6);
\draw  (v7) edge (v8);
\draw  (v7) edge (v9);
\end{tikzpicture}
\caption{Case 2.1.2.}
\label{case 2.1.2.}
\end{minipage}
\end{figure}

\medskip
By the choice of $f_2$ and the symmetry of $2_a$ and $2_b$, the remaining case is:

{\bf Case 2.2:} $f_2(w_1)= f_2(w_2)=2_a$. In particular, this means $w_1w_2\notin E(G)$.
By Lemma~\ref{c1}, $G-w$ has a coloring $g$ satisfying one of the following:\\
(a) $g(w_1)=2_a$ and $g(w_2)\in \{1_a,1_b\}$ or  $g(w_2)=2_a$ and $g(w_1)\in \{1_a,1_b\}$;\\
(b) $\{g(w_1),g(w_2)\}=\{2_a,2_b\}$;\\
(c) $\{g(w_3),g(w_4)\}=\{g(w_5),g(w_6)\}=\{1_a,1_b\}$, and exactly one of $w_1,w_2$ has color $2_a$.

If (a) or (b) occurs, then we again get Case 1. We do not get Case 1 only if (c) occurs and one of $w_1,w_2$ has
$g$-color in $\{3_a,3_b\}$. But then $2_b$ is not present in $B(w,2)$ and we can color $w$ with $2_b$.
\hfill\qed

\section{Cubic graphs}

A {\em good} coloring is a $(1,1,2,2,3,3,4)$-coloring with color $4$ used at most once. By Proposition~\ref{extension},
Theorem~\ref{t1} follows from the following fact.

\begin{thm}\label{1122334}
Every connected cubic graph has a good coloring.
\end{thm}

{\bf Proof.}
Let $G$ be a connected cubic graph with $n\geq 2$ vertices. 
 Since $G$ is connected, it has a non-cut vertex $w$ (simply take a leaf vertex of a spanning tree of $G$). Let $N(w)=\{w_1,w_2,w_3\}$. 

{\bf Case 1:} $0 \le |E(G[\{w_1,w_2,w_3\}])| \le 1$. If $|E(G[\{w_1,w_2,w_3\}])|=0$, then let $G'=G-w+w_2w_3$. If $|E(G[\{w_1,w_2,w_3\}])|=1$, then by symmetry we may assume $w_2w_3 \in E(G)$. Let $G'=G-w$. Note that $G'$ is a connected subcubic graph with vertex $w_1$ of degree at most two. By Theorem~\ref{t2},
  $G'$ has a feasible coloring. Hence by Lemma~\ref{c2},  $G'$ has a feasible coloring $f$ with
 \begin{equation}\label{eeq0} 
   f(w_1)\notin \{3_a,3_b\}.
 \end{equation}    
  Let $N_{G'}(w_1)=\{w_4,w_5\}$, $N_{G'}(w_2)=\{w_3,w_6,w_7\}$, and $N_{G'}(w_3)=\{w_2,w_8,w_9\}.$ It is possible that $|\{w_4,w_5,w_6,w_7,w_8,w_9\}|<6$, but this will not affect the proof below. 
  
 For $j\in\{1,2,3\}$ and $x,y\in V(G)-w$, a $(j,x,y)$-{\em conflict} in $(G,f)$ is the situation that 
  $f(x)=f(y)\in \{j_a,j_b\}$ and $d_G(x,y) \le j$. If $(G,f)$ has no $(j,x,y)$-conflicts for any $j\in\{1,2,3\}$ and $x,y\in V(G)-w$,  
 then we can extend $f$ to a good coloring of $G$ by letting $f(w)=4$. 
 
Suppose now that  $(G,f)$ has a $(j,x,y)$-conflict for some $j\in\{1,2,3\}$ and $x,y\in V(G)-w$ (there could be more than one conflict).
 Then
 \begin{equation}\label{eeq1}
\mbox{\em  $d_{G}(x,y)\leq j<d_{G'}(x,y)$. This means $\{x,y\}\cap \{w_1,w_2,w_3\}\neq \emptyset$ and $j\geq 2$. }
 \end{equation}
 
 Since $w_2w_3\in E(G')$,~(\ref{eeq1}) yields that in each $(j,x,y)$-conflict, one of
  $x$ and $y$ is in $ \{w_1,w_4,w_5\}$ and the other is in $ \{w_2,w_3,w_6,w_7,w_8,w_9\}$. 
  By~(\ref{eeq0}), we have the following two cases.

{\bf Case 1.1:} $f(w_1)\in \{1_a,1_b\}$, say $f(w_1)=1_a$. Then each conflict is a $(3,x,y)$-conflict.

{\bf Case 1.1.1:} There is only one conflict. We may assume it is a $(3,w_4,w_2)$-conflict,
where $f(w_4)=f(w_2)=3_a$. If $f(N_G(w_2)-w)\neq \{1_a,1_b\}$, then we can recolor $w_2$ with one of $1_a$ and $1_b$ and
eliminate the conflict.
If $f(w_3)\neq 1_b$, then we can recolor $w_4$ with $4$ and color $w$ with $1_b$.
So we may assume 
 \begin{equation}\label{eeq11}
\mbox{\em  $f(N_G(w_2)-w)= \{1_a,1_b\}$ \quad and \quad $f(w_3)=1_b$. }
 \end{equation}
Furthermore, if $f(w_5)\neq 1_b$ or $1_a \notin f(N_G(w_3)-w)$, then we can recolor $w_1$ and $w_3$
with the same color $\alpha\in \{1_a,1_b\}$, recolor $w_4$ with $4$ and color $w$ with $\beta\in \{1_a,1_b\}-\alpha$.
Otherwise, some $\gamma\in \{2_a,2_b\}$ is not present on $N(w_3)\cup \{w_5\}$, and by~(\ref{eeq11}) we can
recolor $w_4$ with $4$ and color $w$ with $\gamma$ (See Figure~\ref{case 1.1.1.}).

{\bf Case 1.1.2:} There are two conflicts. By the case and symmetry, we may assume
 $f(w_4)=f(w_2)=3_a$ and $f(w_5)=f(w_3)=3_b$. Applying Lemma~\ref{c2} to vertex $w_2$ and coloring $f$ of
 $G-w$, we obtain a feasible coloring $g$ of $G-w$ such that $g(w_2)=\gamma \notin \{3_a,3_b\}$ and at most one of 
 $w_3,w_4,w_5$ changed its color. 
 
{\bf Case 1.1.2.1:} Neither $w_4$ nor $w_5$ changed its color. Then we color $w_3$ with color $4$, $w$ with a color $\beta \in \{1_a,1_b\}-\gamma$, $w_1$ with a color $\alpha \in \{1_a,1_b\}-\beta$, and use $g$ on other vertices.

{\bf Case 1.1.2.2:} One vertex of $\{w_4,w_5\}$ changed its color. We prove the case when $w_4$ changed its color, say $g(w_4)=\beta \in \{1_a,1_b\}$, the case $w_5$ changed its color is similar. We may assume that 
\begin{equation}\label{gammabeta}
g(w_2)=\gamma \in \{1_a,1_b\}\quad  \mbox{and}\quad  \gamma = \beta,
\end{equation}
since otherwise we color $w_1$ with a color $\alpha \in \{1_a,1_b\}-\beta$, $w$ with a color $\mu \in \{1_a,1_b\}-\alpha$, $w_3$ with color $4$, and use $g$ on other vertices. We may also assume that some vertex, say $w_6 \in N(w_2)-w$, have color $\delta \in \{1_a,1_b\}-\gamma$, since otherwise we recolor $w_2$ with $\delta$ and it contradicts ~\eqref{gammabeta}. We may also assume that $g(\{w_8,w_9\})=\{1_a,1_b\}$, since otherwise we color $w_3$ with a color $\mu \in \{1_a,1_b\}-g(\{w_8,w_9\})$, $w$ with color $4$, and use $f$ on other vertices (See Figure~\ref{case 1.1.2.2.}). Note that $|g(N(w) \cup N(N(w))) \cap \{2_a,2_b\}| \le 1$. Then we color $w_1$ with a color $\alpha \in \{1_a,1_b\}-\beta$, $w_3$ with color $4$, $w$ with a color $\lambda \in \{2_a,2_b\}-g(N(w) \cup N(N(w)))$, and use $g$ on other vertices to obtain a good coloring.

\begin{figure}[ht]\label{f9}
\hspace{5mm}
\begin{minipage}[b]{0.4\textwidth}
\begin{tikzpicture}[scale=0.5, transform shape]

\node[circle, draw=white!0,  inner sep=0pt, minimum size=25pt, font=\huge] (11) at (0,1) {$\gamma$};
\node[circle, draw=white!0,  inner sep=0pt, minimum size=25pt, font=\huge] (12) at (-5,-1) {$1_a$};
\node[circle, draw=white!0,  inner sep=0pt, minimum size=25pt, font=\huge] (14) at (5,-1) {$1_b$};

\node[circle, draw=white!0,  inner sep=0pt, minimum size=25pt, font=\huge] (15) at (6.5,-3) {$2_a /2_b$};

\node[circle, draw=white!0,  inner sep=0pt, minimum size=25pt, font=\huge] (19) at (-1,-3) {$1_a$};

\node[circle, draw=white!0,  inner sep=0pt, minimum size=25pt, font=\huge] (20) at (-6,-3) {$3_a$};

\node[circle, draw=white!0, inner sep=0pt, minimum size=25pt, font=\huge] (24) at (4,-3) {$1_a$};

\node[circle, draw=white!0, inner sep=0pt, minimum size=25pt, font=\huge] (26) at (1,-3) {$1_b$};

\node[circle, draw=white!0, inner sep=0pt, minimum size=25pt, font=\huge] (37) at (0.5,-1) {$3_a$};

\node[circle, draw=white!0, inner sep=0pt, minimum size=25pt, font=\huge] (39) at (-4,-3) {$1_b$};

\node[circle, draw=black!80, inner sep=0pt, minimum size=25pt, font=\huge] (1) at (0,0) {$w$};
\node[circle, draw=black!80, inner sep=0pt, minimum size=25pt, font=\huge] (3) at (-5,-2) { $w_1$};
\node[circle, draw=black!80, inner sep=0pt, minimum size=25pt, font=\huge] (4) at (5,-2) {$w_3$};
\node[circle, draw=black!80, inner sep=0pt, minimum size=25pt, font=\huge] (5) at (-6,-4) {$w_4$};

\node[circle, draw=black!80, inner sep=0pt, minimum size=25pt, font=\huge] (8) at (4,-4) {$w_8$};

\node[circle, draw=black!80, inner sep=0pt, minimum size=25pt, font=\huge] (29) at (-4,-4) {$w_5$};
\node[circle, draw=black!80, inner sep=0pt, minimum size=25pt, font=\huge] (v7) at (6,-4) {$w_9$};

\draw  (1) edge (3);
\draw  (1) edge (4);
\draw  (4) edge (8);
\draw  (4) edge (v7);
\draw  (29) edge (3);
\draw  (5) edge (3);

\node (v1) at (-6.5,-6) {};
\node (v2) at (-5.5,-6) {};
\node (v3) at (-4.5,-6) {};
\node (v4) at (-3.5,-6) {};
\node (v5) at (3.5,-6) {};
\node (v6) at (4.5,-6) {};
\node (v8) at (5.5,-6) {};
\node (v9) at (6.5,-6) {};
\draw  (5) edge (v1);
\draw  (5) edge (v2);
\draw  (29) edge (v3);
\draw  (29) edge (v4);
\draw  (8) edge (v5);
\draw  (8) edge (v6);
\draw  (v7) edge (v8);
\draw  (v7) edge (v9);
\node[circle, draw=black!80, inner sep=0pt, minimum size=25pt, font=\huge] (v10) at (0,-2) {$w_2$};
\draw  (1) edge (v10);
\node[circle, draw=black!80, inner sep=0pt, minimum size=25pt, font=\huge] (v11) at (-1,-4) {$w_6$};
\node[circle, draw=black!80, inner sep=0pt, minimum size=25pt, font=\huge] (v12) at (1,-4) {$w_7$};
\draw  (v10) edge (v11);
\draw  (v10) edge (v12);
\draw  (v10) edge (4);

\node (v13) at (-1.5,-6) {};
\node (v14) at (-0.5,-6) {};
\node (v15) at (0.5,-6) {};
\node (v16) at (1.5,-6) {};
\draw  (v11) edge (v13);
\draw  (v11) edge (v14);
\draw  (v12) edge (v15);
\draw  (v12) edge (v16);

\end{tikzpicture}

\caption{Case 1.1.1.}
\label{case 1.1.1.}
\end{minipage}
\hspace{15mm}
\begin{minipage}[b]{0.4\textwidth}
\begin{tikzpicture}[scale=0.5, transform shape]
\node[circle, draw=white!0,  inner sep=0pt, minimum size=25pt, font=\huge] (11) at (0,1) {$\lambda$};
\node[circle, draw=white!0,  inner sep=0pt, minimum size=25pt, font=\huge] (12) at (-5,-1) {$1_a$};
\node[circle, draw=white!0,  inner sep=0pt, minimum size=25pt, font=\huge] (14) at (5,-1) {$3_b$};

\node[circle, draw=white!0,  inner sep=0pt, minimum size=25pt, font=\huge] (15) at (6,-3) {$1_b$};

\node[circle, draw=white!0,  inner sep=0pt, minimum size=25pt, font=\huge] (19) at (-1,-3) {$\delta$};

\node[circle, draw=white!0,  inner sep=0pt, minimum size=25pt, font=\huge] (20) at (-6,-3) {$\beta$};

\node[circle, draw=white!0, inner sep=0pt, minimum size=25pt, font=\huge] (24) at (4,-3) {$1_a$};

\node[circle, draw=white!0, inner sep=0pt, minimum size=25pt, font=\huge] (26) at (1.5,-3) {$2_a /2_b$};

\node[circle, draw=white!0, inner sep=0pt, minimum size=25pt, font=\huge] (37) at (0.5,-1) {$\gamma$};

\node[circle, draw=white!0, inner sep=0pt, minimum size=25pt, font=\huge] (39) at (-4,-3) {$3_b$};

\node[circle, draw=black!80, inner sep=0pt, minimum size=25pt, font=\huge] (1) at (0,0) {$w$};
\node[circle, draw=black!80, inner sep=0pt, minimum size=25pt, font=\huge] (3) at (-5,-2) { $w_1$};
\node[circle, draw=black!80, inner sep=0pt, minimum size=25pt, font=\huge] (4) at (5,-2) {$w_3$};
\node[circle, draw=black!80, inner sep=0pt, minimum size=25pt, font=\huge] (5) at (-6,-4) {$w_4$};

\node[circle, draw=black!80, inner sep=0pt, minimum size=25pt, font=\huge] (8) at (4,-4) {$w_8$};

\node[circle, draw=black!80, inner sep=0pt, minimum size=25pt, font=\huge] (29) at (-4,-4) {$w_5$};
\node[circle, draw=black!80, inner sep=0pt, minimum size=25pt, font=\huge] (v7) at (6,-4) {$w_9$};

\draw  (1) edge (3);
\draw  (1) edge (4);
\draw  (4) edge (8);
\draw  (4) edge (v7);
\draw  (29) edge (3);
\draw  (5) edge (3);

\node (v1) at (-6.5,-6) {};
\node (v2) at (-5.5,-6) {};
\node (v3) at (-4.5,-6) {};
\node (v4) at (-3.5,-6) {};
\node (v5) at (3.5,-6) {};
\node (v6) at (4.5,-6) {};
\node (v8) at (5.5,-6) {};
\node (v9) at (6.5,-6) {};
\draw  (5) edge (v1);
\draw  (5) edge (v2);
\draw  (29) edge (v3);
\draw  (29) edge (v4);
\draw  (8) edge (v5);
\draw  (8) edge (v6);
\draw  (v7) edge (v8);
\draw  (v7) edge (v9);
\node[circle, draw=black!80, inner sep=0pt, minimum size=25pt, font=\huge] (v10) at (0,-2) {$w_2$};
\draw  (1) edge (v10);
\node[circle, draw=black!80, inner sep=0pt, minimum size=25pt, font=\huge] (v11) at (-1,-4) {$w_6$};
\node[circle, draw=black!80, inner sep=0pt, minimum size=25pt, font=\huge] (v12) at (1,-4) {$w_7$};
\draw  (v10) edge (v11);
\draw  (v10) edge (v12);
\draw  (v10) edge (4);

\node (v13) at (-1.5,-6) {};
\node (v14) at (-0.5,-6) {};
\node (v15) at (0.5,-6) {};
\node (v16) at (1.5,-6) {};
\draw  (v11) edge (v13);
\draw  (v11) edge (v14);
\draw  (v12) edge (v15);
\draw  (v12) edge (v16);
\end{tikzpicture}
\caption{Case 1.1.2.2.}
\label{case 1.1.2.2.}
\end{minipage}
\end{figure}

{\bf Case 1.2:} $f(w_1)\in \{2_a,2_b\}$, say $f(w_1)=2_a$. Since we cannot switch to Case 1.1, we need
$\{f(w_4),f(w_5)\}= \{1_a,1_b\}$. So the only possible conflict is a $(2,w_1,y)$-conflict, where $y\in \{w_2,w_3\}$.
We may assume $f(w_2)=2_a$. Then we recolor $w_1$ with $4$ and color $w$ with
$\alpha\in  \{1_a,1_b\}-f(w_3)$.

{\bf Case 2:} $|E(G[\{w_1,w_2,w_3\}])|=2$, say $w_1w_2 \in E(G)$ and $w_2w_3 \in E(G)$. We obtain a good coloring $g$ of $G$ by using $f$ on $G-w$ and assigning color $4$ to $w$. Note that adding $w$ back will not create conflicts because the distance between any two vertices in $G-w$ remains the same.

{\bf Case 3:} $G[\{w_1,w_2,w_3\}]=K_3$. Then $G=K_4$, and $K_4$ has a good coloring. \hfill\qed







\section{Proof of  Lemma~\ref{c2} }
Recall the claim of the lemma:

\medskip\noindent
{\bf Lemma~\ref{c2}.} {\em
Let $G$ be a subcubic graph and $f$ be a feasible coloring of $G$. Suppose there is a $2$-vertex $u \in V(G)$ with $N(u)=\{u_1,u_2\}$. If $f(u) \in \{3_a,3_b\}$, then we can recolor some vertices of $G$ so that the resulting coloring $g$ is feasible and satisfies the following:\\
(a) $g(u) \notin \{3_a,3_b\}$, and\\
(b) at most one vertex is recolored into $3_a$ or $3_b$, and this vertex (if there is such a vertex) is
at distance at most $3$ from $u$ and has degree $3$ in $G$, and at most one vertex of $f$-color $3_a$ or $3_b$ apart from $u$ is recolored into some other color, and this vertex (if there is such a vertex) has new color in $\{1_a, 1_b\}$. }

\bigskip
{\bf Proof.} Without loss of generality, we assume that $f(u)=3_a$. If $\{f(u_1),f(u_2)\} \neq \{1_a,1_b\}$, then we recolor $u$ with a color $x \in \{1_a,1_b\}-\{f(u_1),f(u_2)\}$ to obtain a coloring satisfying $(a)\text{ and }(b)$. Thus we may assume 
\begin{equation}\label{1a1b}
f(u_1)=1_a \quad  \mbox{and}\quad   f(u_2)=1_b.
\end{equation}
Let $G_1$ denote the subgraph of $G$ induced by the vertices of colors $1_a$ and $1_b$. If $u_1$ and $u_2$ are in distinct
components of $G_1$, then after switching the colors in the component of $G_1$ containing $u_2$, we obtain a coloring
contradicting~\eqref{1a1b}. Thus we may assume
\begin{equation}\label{1a1b'}
  \mbox{\em $G$ has a $1_a,1_b$-colored $u_1,u_2$-path $P_u$.}
\end{equation}

{\bf Case 1:} $u_1u_2\in E(G)$. If $|N(u_1)|=3$, then let $u_3\in N(u_1)-\{u,u_2\}$. Similarly, if
$|N(u_2)|=3$, then let $u_4\in N(u_2)-\{u,u_1\}$. If $\{2_a, 2_b\} \nsubseteq f(N(u_1)\cup N(u_2))$, then after recoloring $u$ with a color $x \in \{2_a, 2_b\}-f(N(u_1)\cup N(u_2))$ we obtain a
coloring satisfying $(a)\text{ and }(b)$. By symmetry, we may assume 
\begin{equation}\label{u3u41}
  \mbox{\em
$|N(u_1)|=|N(u_2)|=3$,  \quad $f(u_3)=2_a$ \quad and \quad $f(u_4) = 2_b$.}
\end{equation}
If $1_b \notin f(N(u_3))$, then we can recolor $u_3$ with $1_b$ and $u$ with $2_a$ to obtain a
coloring satisfying $(a)\text{ and }(b)$. So we may assume $1_b \in f(N(u_3))$. Similarly, we may assume $1_a \in f(N(u_4))$. If $|N(u_3)|=2$ or $1_a \notin f(N(u_3)-\{u_1\})$, then we can recolor $u_3$ with $1_a$, $u_1$ with $2_a$, and $u$ with $1_a$ to obtain a coloring satisfying $(a)\text{ and }(b)$. So we may assume
\begin{equation}\label{u5u6}
|N(u_3)|=3 \text{ and let }u_5,u_6 \in N(u_3)-\{u_1\} \text{ with } f(u_5)=1_a, f(u_6)=1_b.
\end{equation}
  Similarly, we may assume 
\begin{equation}\label{u7u8}
|N(u_4)|=3 \text{ and let } u_7,u_8 \in N(u_4)-\{u_2\} \text{ with } f(u_7)=1_a, f(u_8)=1_b.
\end{equation}  
  
{\bf Case 1.1:} $u_5=u_7 \text{ and } u_6=u_8$. If $1_b \notin f(N(u_5))$, then we can recolor $u_5$ with $1_b$, $u_3$ with $1_a$, $u_1$ with $2_a$, and $u$ with $1_a$ to obtain a coloring satisfying $(a)\text{ and }(b)$. So we may assume $1_b \in f(N(u_5))$. Similarly, we may assume $1_a \in f(N(u_6))$. Then we can recolor $u_1$ with $3_a$ and $u$ with $1_a$ to obtain a coloring satisfying $(a)\text{ and }(b)$.

{\bf Case 1.2:} $u_5=u_7$ or $u_6=u_8$, but not both. By symmetry, we may assume $u_6 = u_8$ and $u_5 \neq u_7$. It is possible that $u_5u_6 \in E(G)$ or $u_6u_7 \in E(G)$, but this will not affect the proof below. 

Similarly to Case 1.1, we may assume 
\begin{equation}\label{u5u6u7}
1_b \in f(N(u_5)), \quad 1_a \in f(N(u_6)) \quad  \mbox{and}\quad  1_b \in f(N(u_7)).
\end{equation}
Since $3_a \notin f(N(u_6))$, we can also assume $3_a \in f(N(u_5))$, because otherwise we recolor $u_1$ with $3_a$ and $u$ with $1_a$ to obtain a coloring satisfying $(a)\text{ and }(b)$. With \eqref{u3u41} and \eqref{u5u6u7}, we have $f(N(u_5))=\{1_b,2_a,3_a\}$. However, we can recolor $u_1$ with $3_b$ and $u$ with $1_a$ to obtain a coloring satisfying $(a)\text{ and }(b)$.


{\bf Case 1.3:} $u_5 \neq u_7$ and $u_6 \neq u_8$. Then $N(u_3) \cap N(u_4) = \emptyset$ and $d(u_3,u_4)\ge 3$. Similarly to Case 1.2, $\{1_a,1_b,3_a,3_b\} \subseteq f(N(u_5) \cup N(u_6) - \{u_3\})$ (See Figure~\ref{2-case 1.3.}). Therefore, we can recolor $u_3$ with $2_b$ and $u$ with $2_a$ to obtain a coloring satisfying $(a)\text{ and }(b)$.

\begin{figure}[ht]\label{f2}
\centering
\begin{minipage}[b]{0.3\textwidth}
\begin{tikzpicture}[scale=0.5, transform shape]

\node[circle, draw=white!0,  inner sep=0pt, minimum size=25pt, font=\huge] (11) at (1.1,0.7) {$3_a$};
\node[circle, draw=white!0,  inner sep=0pt, minimum size=25pt, font=\huge] (12) at (-2,-1) {$1_a$};
\node[circle, draw=white!0,  inner sep=0pt, minimum size=25pt, font=\huge] (14) at (3,-1) {$1_b$};

\node[circle, draw=white!0,  inner sep=0pt, minimum size=25pt, font=\huge] (17) at (-2.85,-3.5) {$2_a$};
\node[circle, draw=white!0,  inner sep=0pt, minimum size=25pt, font=\huge] (18) at (-0.45,-5.1) {$1_b$};
\node[circle, draw=white!0,  inner sep=0pt, minimum size=25pt, font=\huge] (19) at (-3.5,-5.1) {$1_a$};

\node[circle, draw=white!0,  inner sep=0pt, minimum size=25pt, font=\huge] (22) at (3.8,-3.5) {$2_b$};

\node[circle, draw=white!0,  inner sep=0pt, minimum size=25pt, font=\huge] (23) at (4.5,-5.1) {$1_b$};
\node[circle, draw=white!0,  inner sep=0pt, minimum size=25pt, font=\huge] (24) at (1.8,-5.1) {$1_a$};

\node[circle, draw=white!0,  inner sep=0pt, minimum size=25pt, font=\huge] (25) at (-3.5,-8.5) {$3_a$};
\node[circle, draw=white!0,  inner sep=0pt, minimum size=25pt, font=\huge] (26) at (-2.5,-8.5) {$1_b$};
\node[circle, draw=white!0,  inner sep=0pt, minimum size=25pt, font=\huge] (27) at (-1.5,-8.5) {$1_a$};
\node[circle, draw=white!0,  inner sep=0pt, minimum size=25pt, font=\huge] (28) at (-0.5,-8.5) {$3_b$};

\node[circle, draw=black!80, inner sep=0pt, minimum size=25pt, font=\huge] (1) at (0.5,0) {$u$};
\node[circle, draw=black!80, inner sep=0pt, minimum size=25pt, font=\huge] (3) at (-2,-2) { $u_1$};
\node[circle, draw=black!80, inner sep=0pt, minimum size=25pt, font=\huge] (4) at (3,-2) {$u_2$};
\node[circle, draw=black!80, inner sep=0pt, minimum size=25pt, font=\huge] (5) at (-2,-4) {$u_3$};
\node[circle, draw=black!80, inner sep=0pt, minimum size=25pt, font=\huge] (6) at (-3,-6) {$u_5$};
\node[circle, draw=black!80, inner sep=0pt, minimum size=25pt, font=\huge] (7) at (-1,-6) {$u_6$};
\node[circle, draw=black!80, inner sep=0pt, minimum size=25pt, font=\huge] (8) at (3,-4) {$u_4$};
\draw  (1) edge (3);
\draw  (1) edge (4);
\draw  (3) edge (4);
\draw  (3) edge (5);
\draw  (4) edge (8);
\draw  (5) edge (6);
\draw  (5) edge (7);

\node[circle, draw=black!80,  inner sep=0pt, minimum size=25pt, font=\huge] (v1) at (2,-6) {$u_7$};
\node[circle, draw=black!80,  inner sep=0pt, minimum size=25pt, font=\huge] (v2) at (4,-6) {$u_8$};

\draw  (8) edge (v1);
\draw  (8) edge (v2);
\node (v3) at (-3.5,-8) {};
\node (v4) at (-2.5,-8) {};
\node (v5) at (-1.5,-8) {};
\node (v6) at (-0.5,-8) {};
\draw  (6) edge (v3);
\draw  (6) edge (v4);
\draw  (7) edge (v5);
\draw  (7) edge (v6);

\end{tikzpicture}

\caption{Case 1.3.}
\label{2-case 1.3.}
\end{minipage}
\hspace{15mm}
\begin{minipage}[b]{0.4\textwidth}
\begin{tikzpicture}[scale=0.4, transform shape]
\node[circle, draw=white!0,  inner sep=0pt, minimum size=25pt, font=\huge] (11) at (0,1) {$3_a$};
\node[circle, draw=white!0,  inner sep=0pt, minimum size=25pt, font=\huge] (12) at (-4.15,-1) {$1_a$};
\node[circle, draw=white!0,  inner sep=0pt, minimum size=25pt, font=\huge] (14) at (4,-1) {$1_b$};


\node[circle, draw=white!0,  inner sep=0pt, minimum size=25pt, font=\huge] (17) at (-6,-3) {$2_a$};
\node[circle, draw=white!0,  inner sep=0pt, minimum size=25pt, font=\huge] (18) at (-2,-3) {$1_b$};
\node[circle, draw=white!0,  inner sep=0pt, minimum size=25pt, font=\huge] (19) at (-3.05,-4.85) {$1_a$};

\node[circle, draw=white!0,  inner sep=0pt, minimum size=25pt, font=\huge] (20) at (-1,-4.85) {$2_a$};

\node[circle, draw=white!0, inner sep=0pt, minimum size=25pt, font=\huge] (21) at (5.5,-3) {$2_b$};
\node[circle, draw=white!0, inner sep=0pt, minimum size=25pt, font=\huge] (22) at (-7.05,-4.85) {$2_b$};

\node[circle, draw=white!0, inner sep=0pt, minimum size=25pt, font=\huge] (23) at (-5,-4.85) {$1_b$};
\node[circle, draw=white!0, inner sep=0pt, minimum size=25pt, font=\huge] (24) at (2.5,-3) {$1_a$};

\node[circle, draw=white!0, inner sep=0pt, minimum size=25pt, font=\huge] (25) at (-4.5,-8.5) {$3_a$};
\node[circle, draw=white!0, inner sep=0pt, minimum size=25pt, font=\huge] (26) at (-6.5,-8.5) {$1_b$};
\node[circle, draw=white!0, inner sep=0pt, minimum size=25pt, font=\huge] (27) at (-7.5,-8.5) {$1_a$};
\node[circle, draw=white!0, inner sep=0pt, minimum size=25pt, font=\huge] (28) at (-0.5,-8.5) {$3_b$};
\node[circle, draw=white!0, inner sep=0pt, minimum size=25pt, font=\huge] (30) at (-5.5,-8.5) {$1_a$};
\node[circle, draw=white!0, inner sep=0pt, minimum size=25pt, font=\huge] (31) at (-2.5,-8.5) {$1_b$};
\node[circle, draw=white!0, inner sep=0pt, minimum size=25pt, font=\huge] (32) at (-1.5,-8.5) {$1_a$};
\node[circle, draw=white!0, inner sep=0pt, minimum size=15pt, font=\huge] (33) at (-3.5,-8.5) {$2_b$};
\node[circle, draw=white!0, inner sep=0pt, minimum size=15pt, font=\huge] (34) at (-0.5,-9.5) {$2_b$};
\node[circle, draw=white!0, inner sep=0pt, minimum size=25pt, font=\huge] (35) at (-3.5,-9.5) {$3_b$};

\node[circle, draw=black!80, inner sep=0pt, minimum size=25pt, font=\huge] (1) at (0,0) {$u$};
\node[circle, draw=black!80, inner sep=0pt, minimum size=25pt, font=\huge] (3) at (-4,-2) { $u_1$};
\node[circle, draw=black!80, inner sep=0pt, minimum size=25pt, font=\huge] (4) at (4,-2) {$u_2$};
\node[circle, draw=black!80, inner sep=0pt, minimum size=25pt, font=\huge] (5) at (-6,-4) {$u_3$};
\node[circle, draw=black!80, inner sep=0pt, minimum size=36pt, font=\huge] (6) at (-7,-6) {$u_7$};
\node[circle, draw=black!80, inner sep=0pt, minimum size=36pt, font=\huge] (7) at (-5,-6) {$u_8$};
\node[circle, draw=black!80, inner sep=0pt, minimum size=25pt, font=\huge] (8) at (2.5,-4) {$u_5$};
\draw  (1) edge (3);
\draw  (1) edge (4);

\draw  (3) edge (5);
\draw  (4) edge (8);
\draw  (5) edge (6);
\draw  (5) edge (7);

\node (v3) at (-7.5,-8) {};
\node (v4) at (-6.5,-8) {};
\node (v5) at (-5.5,-8) {};
\node (v6) at (-4.5,-8) {};
\draw  (6) edge (v3);
\draw  (6) edge (v4);
\draw  (7) edge (v5);
\draw  (7) edge (v6);
\node[circle, draw=black!80, inner sep=0pt, minimum size=25pt, font=\huge] (29) at (-2,-4) {$u_4$};
\node[circle, draw=black!80, inner sep=0pt, minimum size=25pt, font=\huge] (v7) at (5.5,-4) {$u_6$};
\draw  (3) edge (29);
\draw  (4) edge (v7);
\node[circle, draw=black!80, inner sep=0pt, minimum size=36pt, font=\huge] (v1) at (-3,-6) {$u_9$};
\node[circle, draw=black!80, inner sep=0pt, minimum size=30pt, font=\huge] (v2) at (-1,-6) {$u_{10}$};
\draw  (29) edge (v1);
\draw  (29) edge (v2);
\node (v8) at (-3.5,-8) {};
\node (v9) at (-2.5,-8) {};
\node (v10) at (-1.5,-8) {};
\node (v11) at (-0.5,-8) {};
\draw  (v1) edge (v8);
\draw  (v1) edge (v9);
\draw  (v2) edge (v10);
\draw  (v2) edge (v11);
\end{tikzpicture}
\caption{Case 2.1.}
\label{2-case 2.1.}
\end{minipage}
\end{figure}


{\bf Case 2:} $u_1u_2\notin E(G)$. If $\{2_a, 2_b\} \nsubseteq f(N(u_1)\cup N(u_2))$, then after recoloring $u$ with a color $x \in \{2_a, 2_b\}-f(N(u_1)\cup N(u_2))$ we obtain a coloring satisfying $(a)\text{ and }(b)$. With \eqref{1a1b'}, we may assume that 
\begin{equation}\label{u1u3u4}
N(u_1)=\{u,u_3,u_4\}, \quad f(u_3) = 2_a, \quad f(u_4) = 1_b,
\end{equation}
\begin{equation}\label{u2u5u6}
N(u_2)=\{u,u_5,u_6\}, \quad f(u_5) = 1_a \quad  \mbox{and}\quad  f(u_6) = 2_b.
\end{equation}

If $u_3u_4 \in E(G)$, then $1_a \in f(N(u_4)-\{u_1,u_3\})$ because of \eqref{1a1b'}. We also have $2_b \in f(N(u_3)-\{u_1,u_4\})$ because otherwise we can recolor $u_1$ with $2_b$ and $u$ with $1_a$ to obtain a coloring satisfying $(a)\text{ and }(b)$. Thus, we may assume $|N(u_3)|=|N(u_4)|=3$ and let $u_7 \in N(u_3)-\{u_1,u_4\}, u_8 \in N(u_4)-\{u_1,u_3\}$, $f(u_7)=2_b$, and $f(u_8) = 1_a$. Then, we can recolor $u_1$ with $2_a$, $u_3$ with $1_a$, and $u$ with $1_a$ to obtain a coloring satisfying $(a)\text{ and }(b)$. Because of symmetry, we may assume 
\begin{equation}\label{u3u4u5u6}
u_3u_4 \notin E(G)\quad  \mbox{and}\quad  u_5u_6 \notin E(G).
\end{equation}
If $1_b \notin f(N(u_3))$, then we recolor $u_3$ with $1_b$ and $u$ with $2_a$ to obtain a coloring satisfying $(a)\text{ and }(b)$. With \eqref{1a1b'}, we may assume that
\begin{equation}\label{1bu31au4}
1_b \in f(N(u_3))\quad  \mbox{and}\quad  1_a \in f(N(u_4)).
\end{equation}
If $2_b \notin f(B(u_1,2))$, then we can recolor $u_1$ with $2_b$ and $u$ with $1_a$ to obtain a coloring satisfying $(a)\text{ and }(b)$. Thus, we may assume 
\begin{equation}\label{2bu3u4}
2_b \in f(N(u_3)) \cup f(N(u_4)). 
\end{equation}
If $1_a \notin f(N(u_3)-\{u_1\})$ and $2_a \notin f(N(u_4))$, then we can recolor $u_3$ with $1_a$, $u_1$ with $2_a$, and $u$ with $1_a$ to obtain a coloring satisfying $(a)\text{ and }(b)$. Thus, we may assume 
\begin{equation}\label{u3u43}
|N(u_3)|=|N(u_4)|=3
\end{equation}
and 
\begin{equation}\label{u3u4}
1_a \in f(N(u_3)-\{u_1\})\text{ or }2_a \in f(N(u_4)).
\end{equation}
Let $\{u_7,u_8\} \in N(u_3)$, $\{u_9,u_{10}\} \in N(u_4)$. By \eqref{1bu31au4}, we may assume
\begin{equation}\label{u8u9}
f(u_8)=1_b \quad  \mbox{and}\quad  f(u_9)=1_a. 
\end{equation}
By \eqref{2bu3u4} and \eqref{u3u4}, we have
\begin{equation}\label{cases}
\mbox{\em either 
$f(u_7)=2_b$ and $f(u_{10})=2_a$  or  $f(u_7)=1_a$  and  $f(u_{10})=2_b$.}
\end{equation}
If $3_a \notin f(B(u_1,3)-\{u\})$, then we can recolor $u_1$ with $3_a$ and $u$ with $1_a$ to obtain a coloring satisfying $(a)\text{ and }(b)$. Thus, we may assume
\begin{equation}\label{3a}
3_a \in f(B(u_1,3)-\{u\}).
\end{equation}
Similarly, we may assume
\begin{equation}\label{3b}
3_b \in f(B(u_1,3)-\{u\}).
\end{equation}

{\bf Case 2.1:} $f(u_7)=2_b$ and $f(u_{10})=2_a$. By \eqref{u3u4u5u6} and $|N(u_2)|=3$, we have $$\{u_8,u_{10}\} \cap (\{u_i: i \in [6]\} \cup \{u\}) = \emptyset.$$ It is possible that $u_9=u_5$ or $u_7=u_6$, but this will not affect the proof below. 

If $2_b \notin f(B(u_4,2))$, then we can recolor $u_4$ with $2_b$, $u_1$ with $1_b$, and $u$ with $1_a$ to obtain a coloring satisfying $(a)\text{ and }(b)$. Thus, we may assume 
\begin{equation}\label{2bnnu4}
2_b \in f(B(u_4,2)).
\end{equation}
If $1_a \notin f(N(u_{10}))$, then we can recolor $u_{10}$ with $1_a$ and it contradicts \eqref{u3u4}. Thus, we may assume 
\begin{equation}\label{1anu10}
1_a \in f(N(u_{10})). 
\end{equation}
We may also assume 
\begin{equation}\label{1a1bu7}
f(N(u_7)-\{u_3\}) = \{1_a,1_b\},
\end{equation}
because otherwise we can recolor $u_7$ with a color $x \in \{1_a,1_b\}-f(N(u_7)-\{u_1\})$ and it contradicts \eqref{cases}.
By \eqref{3a} and \eqref{3b}, we know that 
\begin{equation}\label{3a3b78910}
\{3_a,3_b\} \subseteq f(N(u_7) \cup N(u_8) \cup N(u_9) \cup N(u_{10})).
\end{equation}
If $\{3_a,3_b\} \subseteq f(N(u_7) \cup N(u_8))$, then by \eqref{1a1bu7} we have $f(N(u_8))=\{2_a,3_a,3_b\}$. Then, we can recolor $u_8$ with $1_a$, $u_3$ with $1_b$, and $u$ with $2_a$ to obtain a coloring satisfying $(a)\text{ and }(b)$. By symmetry, we may assume 
\begin{equation}\label{3bu7u8}
3_b \notin f(N(u_7) \cup N(u_8)).
\end{equation}
 By \eqref{3a3b78910} and \eqref{3bu7u8}, we know that $3_b \in f(N(u_9) \cup N(u_{10}))$. By \eqref{1a1b'}, $1_b \in f(N(u_9)-\{u_4\})$. With \eqref{2bnnu4}, \eqref{1anu10}, and $2_b \notin f(\{u,u_1,u_3,u_9,u_{10}\})$ we know that $$f(N(u_9) \cup N(u_{10})-\{u_4\}) = \{1_a,1_b,2_b,3_b\}\text{, hence }1_b \notin f(N(u_{10})-\{u_4\}) \text{ (See Figure~\ref{2-case 2.1.})}.$$ Therefore, we can recolor $u_{10}$ with $1_b$, $u_4$ with $2_a$, $u_3$ with $1_a$, $u_1$ with $1_b$, and $u$ with $1_a$ to obtain a coloring satisfying $(a)\text{ and }(b)$.


{\bf Case 2.2:} $f(u_7)=1_a$ and $f(u_{10})=2_b$. If $1_a \notin f(N(u_6))$, then we can recolor $u_6$ with $1_a$ and $u$ with $2_b$ to obtain a coloring satisfying $(a)\text{ and }(b)$. Thus, we may assume 
\begin{equation}\label{1au6}
1_a \in f(N(u_6)-\{u_2\}).
\end{equation}

Since some $u_i$ and $u_j$ may coincide, several cases are considered below.

{\bf Case 2.2.1}: $u_3u_5 \in E(G)$, i.e., $u_7=u_5$. It is possible that $u_4u_6 \in E(G)$, or $u_4u_5 \in E(G)$, or $\{u_4u_5, u_4u_6\} \subseteq E(G)$, but this will not affect the proof below. By \eqref{1a1b'}, 
\begin{equation}\label{1bu9}
1_b \in f(N(u_9)-\{u_4\}),
\end{equation}
and
\begin{equation}\label{1bu5}
1_b \in f(N(u_5)-\{u_2\}).
\end{equation}
If $1_a \notin f(N(u_{10})-\{u_4\})$, then we can recolor $u_{10}$ with $1_a$ and it contradicts \eqref{cases}. Thus, we may assume 
\begin{equation}\label{1au10}
1_a \in f(N(u_{10})-\{u_4\}).
\end{equation}
If $1_a \notin f(N(u_8))$, then we can recolor $u_8$ with $1_a$, $u_3$ with $1_b$, and $u$ with $2_a$ to obtain a coloring satisfying $(a)\text{ and }(b)$. If $2_b \notin f(N(u_8))$, then we can recolor $u_3$ with $2_b$ and $u$ with $2_a$ to obtain a coloring satisfying $(a)\text{ and }(b)$. Thus, we may assume 
\begin{equation}\label{u8}
f(N(u_8))=\{1_a,2_a,2_b\}.
\end{equation}
By \eqref{3a}, \eqref{3b}, \eqref{1bu9}, \eqref{1bu5}, \eqref{1au10}, and \eqref{u8}, we have 
\begin{equation}\label{1a1b3a3b}
\{1_a,1_b,3_a,3_b\} \subseteq f(N(u_{9}) \cup N(u_{10})-\{u_4\}).
\end{equation}
By \eqref{1a1b3a3b}, $1_b \notin f(N(u_{10})-\{u_4\}),$ and $2_b \notin f(B(u_4,2)-\{u_{10}\})$ (See Figure~\ref{2-case 2.2.1.}). Then, we can recolor $u_{10}$ with $1_b$, $u_4$ with $2_b$, $u_1$ with $1_b$, and $u$ with $1_a$ to obtain a coloring satisfying $(a)\text{ and }(b)$. 

With Case 2.2.1 handled, from now on by symmetry we may assume
\begin{equation}\label{u3u5u4u6}
u_3u_5 \notin E(G)\quad  \mbox{and}\quad  u_4u_6 \notin E(G).
\end{equation}

\begin{figure}[ht]\label{f3}
\hspace{5mm}
\begin{minipage}[b]{0.4\textwidth}
\begin{tikzpicture}[scale=0.5, transform shape]
\node[circle, draw=white!0,  inner sep=0pt, minimum size=25pt, font=\huge] (11) at (0,1) {$3_a$};
\node[circle, draw=white!0,  inner sep=0pt, minimum size=25pt, font=\huge] (12) at (-4.15,-1) {$1_a$};
\node[circle, draw=white!0,  inner sep=0pt, minimum size=25pt, font=\huge] (14) at (4,-1) {$1_b$};


\node[circle, draw=white!0,  inner sep=0pt, minimum size=25pt, font=\huge] (18) at (2.5,-6) {$1_b$};

\node[circle, draw=white!0, inner sep=0pt, minimum size=25pt, font=\huge] (21) at (5.5,-3) {$2_b$};
\node[circle, draw=white!0, inner sep=0pt, minimum size=25pt, font=\huge] (22) at (-5,-5) {$2_b$};

\node[circle, draw=white!0, inner sep=0pt, minimum size=25pt, font=\huge] (24) at (2.5,-3) {$1_a$};

\node[circle, draw=white!0, inner sep=0pt, minimum size=25pt, font=\huge] (25) at (-7.5,-8.5) {$3_a$};
\node[circle, draw=white!0, inner sep=0pt, minimum size=25pt, font=\huge] (26) at (-6.5,-8.5) {$1_b$};

\node[circle, draw=white!0, inner sep=0pt, minimum size=25pt, font=\huge] (30) at (-5.5,-8.5) {$1_a$};
\node[circle, draw=white!0, inner sep=0pt, minimum size=25pt, font=\huge] (31) at (-6,-3) {$1_b$};
\node[circle, draw=white!0, inner sep=0pt, minimum size=25pt, font=\huge] (32) at (-2.5,-8.5) {$1_a$};

\node[circle, draw=white!0, inner sep=0pt, minimum size=15pt, font=\huge] (34) at (-1.5,-8.5) {$2_b$};
\node[circle, draw=white!0, inner sep=0pt, minimum size=25pt, font=\huge] (35) at (-4.5,-8.5) {$3_b$};

\node[circle, draw=white!0, inner sep=0pt, minimum size=25pt, font=\huge] (38) at (-7,-5) {$1_a$};
\node[circle, draw=white!0, inner sep=0pt, minimum size=25pt, font=\huge] (39) at (-1.5,-5) {$1_b$};
\node[circle, draw=white!0, inner sep=0pt, minimum size=25pt, font=\huge] (40) at (-2,-3) {$2_a$};

\node[circle, draw=black!80, inner sep=0pt, minimum size=25pt, font=\huge] (1) at (0,0) {$u$};
\node[circle, draw=black!80, inner sep=0pt, minimum size=25pt, font=\huge] (3) at (-4,-2) { $u_1$};
\node[circle, draw=black!80, inner sep=0pt, minimum size=25pt, font=\huge] (4) at (4,-2) {$u_2$};
\node[circle, draw=black!80, inner sep=0pt, minimum size=25pt, font=\huge] (5) at (-6,-4) {$u_4$};
\node[circle, draw=black!80, inner sep=0pt, minimum size=36pt, font=\huge] (6) at (-7,-6) {$u_9$};
\node[circle, draw=black!80, inner sep=0pt, minimum size=36pt, font=\huge] (7) at (-5,-6) {$u_{10}$};
\node[circle, draw=black!80, inner sep=0pt, minimum size=25pt, font=\huge] (8) at (2.5,-4) {$u_5$};
\draw  (1) edge (3);
\draw  (1) edge (4);

\draw  (3) edge (5);
\draw  (4) edge (8);
\draw  (5) edge (6);
\draw  (5) edge (7);

\node (v3) at (-7.5,-8) {};
\node (v4) at (-6.5,-8) {};
\node (v5) at (-5.5,-8) {};
\node (v6) at (-4.5,-8) {};
\draw  (6) edge (v3);
\draw  (6) edge (v4);
\draw  (7) edge (v5);
\draw  (7) edge (v6);
\node[circle, draw=black!80, inner sep=0pt, minimum size=25pt, font=\huge] (29) at (-2,-4) {$u_3$};
\node[circle, draw=black!80, inner sep=0pt, minimum size=25pt, font=\huge] (v7) at (5.5,-4) {$u_6$};
\draw  (3) edge (29);
\draw  (4) edge (v7);

\node[circle, draw=black!80, inner sep=0pt, minimum size=30pt, font=\huge] (v2) at (-2,-6) {$u_{8}$};

\draw  (29) edge (v2);

\node (v10) at (-2.5,-8) {};
\node (v11) at (-1.5,-8) {};

\draw  (v2) edge (v10);
\draw  (v2) edge (v11);

\draw  (29) edge (8);
\node (v12) at (2.5,-5.5) {};
\draw  (8) edge (v12);
\node (v13) at (5,-6) {};
\node (v14) at (6,-6) {};
\draw  (v7) edge (v13);
\draw  (v7) edge (v14);

\end{tikzpicture}

\caption{Case 2.2.1.}
\label{2-case 2.2.1.}
\end{minipage}
\hspace{15mm}
\begin{minipage}[b]{0.4\textwidth}
\begin{tikzpicture}[scale=0.5, transform shape]
\node[circle, draw=white!0,  inner sep=0pt, minimum size=25pt, font=\huge] (11) at (0,1) {$3_a$};
\node[circle, draw=white!0,  inner sep=0pt, minimum size=25pt, font=\huge] (12) at (-4.15,-1) {$1_a$};
\node[circle, draw=white!0,  inner sep=0pt, minimum size=25pt, font=\huge] (14) at (4,-1) {$1_b$};


\node[circle, draw=white!0,  inner sep=0pt, minimum size=25pt, font=\huge] (18) at (2.5,-6) {$2_a$};

\node[circle, draw=white!0, inner sep=0pt, minimum size=25pt, font=\huge] (21) at (5.5,-3) {$2_b$};
\node[circle, draw=white!0, inner sep=0pt, minimum size=25pt, font=\huge] (22) at (-1.5,-5) {$2_b$};

\node[circle, draw=white!0, inner sep=0pt, minimum size=25pt, font=\huge] (24) at (2.5,-3) {$1_a$};

\node[circle, draw=white!0, inner sep=0pt, minimum size=25pt, font=\huge] (25) at (-7.5,-8.5) {$3_a$};
\node[circle, draw=white!0, inner sep=0pt, minimum size=25pt, font=\huge] (26) at (-6.5,-8.5) {$1_b$};

\node[circle, draw=white!0, inner sep=0pt, minimum size=25pt, font=\huge] (30) at (-5.5,-8.5) {$1_a$};
\node[circle, draw=white!0, inner sep=0pt, minimum size=25pt, font=\huge] (31) at (-2,-3) {$1_b$};
\node[circle, draw=white!0, inner sep=0pt, minimum size=25pt, font=\huge] (32) at (-2.5,-8.5) {$1_a$};

\node[circle, draw=white!0, inner sep=0pt, minimum size=15pt, font=\huge] (34) at (-4.5,-8.5) {$2_b$};
\node[circle, draw=white!0, inner sep=0pt, minimum size=25pt, font=\huge] (35) at (-1.5,-8.5) {$3_b$};

\node[circle, draw=white!0, inner sep=0pt, minimum size=25pt, font=\huge] (38) at (-7,-5) {$1_a$};
\node[circle, draw=white!0, inner sep=0pt, minimum size=25pt, font=\huge] (39) at (-5,-5) {$1_b$};
\node[circle, draw=white!0, inner sep=0pt, minimum size=25pt, font=\huge] (40) at (-6,-3) {$2_a$};

\node[circle, draw=white!0, inner sep=0pt, minimum size=25pt, font=\huge] (41) at (6,-6.5) {$1_b$};

\node[circle, draw=white!0, inner sep=0pt, minimum size=25pt, font=\huge] (42) at (5,-6.5) {$1_a$};
\node[circle, draw=black!80, inner sep=0pt, minimum size=25pt, font=\huge] (1) at (0,0) {$u$};
\node[circle, draw=black!80, inner sep=0pt, minimum size=25pt, font=\huge] (3) at (-4,-2) { $u_1$};
\node[circle, draw=black!80, inner sep=0pt, minimum size=25pt, font=\huge] (4) at (4,-2) {$u_2$};
\node[circle, draw=black!80, inner sep=0pt, minimum size=25pt, font=\huge] (5) at (-6,-4) {$u_3$};
\node[circle, draw=black!80, inner sep=0pt, minimum size=36pt, font=\huge] (6) at (-7,-6) {$u_7$};
\node[circle, draw=black!80, inner sep=0pt, minimum size=36pt, font=\huge] (7) at (-5,-6) {$u_{8}$};
\node[circle, draw=black!80, inner sep=0pt, minimum size=25pt, font=\huge] (8) at (2.5,-4) {$u_5$};
\draw  (1) edge (3);
\draw  (1) edge (4);

\draw  (3) edge (5);
\draw  (4) edge (8);
\draw  (5) edge (6);
\draw  (5) edge (7);

\node (v3) at (-7.5,-8) {};
\node (v4) at (-6.5,-8) {};
\node (v5) at (-5.5,-8) {};
\node (v6) at (-4.5,-8) {};
\draw  (6) edge (v3);
\draw  (6) edge (v4);
\draw  (7) edge (v5);
\draw  (7) edge (v6);
\node[circle, draw=black!80, inner sep=0pt, minimum size=25pt, font=\huge] (29) at (-2,-4) {$u_4$};
\node[circle, draw=black!80, inner sep=0pt, minimum size=25pt, font=\huge] (v7) at (5.5,-4) {$u_6$};
\draw  (3) edge (29);
\draw  (4) edge (v7);

\node[circle, draw=black!80, inner sep=0pt, minimum size=30pt, font=\huge] (v2) at (-2,-6) {$u_{10}$};

\draw  (29) edge (v2);

\node (v10) at (-2.5,-8) {};
\node (v11) at (-1.5,-8) {};

\draw  (v2) edge (v10);
\draw  (v2) edge (v11);

\draw  (29) edge (8);
\node (v12) at (2.5,-5.5) {};
\draw  (8) edge (v12);
\node (v13) at (5,-6) {};
\node (v14) at (6,-6) {};
\draw  (v7) edge (v13);
\draw  (v7) edge (v14);

\end{tikzpicture}
\caption{Case 2.2.2.}
\label{2-case 2.2.2.}
\end{minipage}
\end{figure}

{\bf Case 2.2.2:} $\{u_3u_5, u_4u_6\} \cap E(G) = \emptyset$ and $u_4u_5 \in E(G)$, i.e., $u_9=u_5$. If $2_a \notin f(N(u_5) \cup N(u_6))$, then we can recolor $u_2$ with $2_a$ and $u$ with $1_b$ to obtain a coloring satisfying $(a)\text{ and }(b)$. If $1_b \notin f(N(u_6)-\{u_2\})$ and $2_b \notin f(N(u_5)-\{u_2,u_4\})$, then we can recolor $u_6$ with $1_b$, $u_2$ with $2_b$, and $u$ with $1_b$ to obtain a coloring satisfying $(a)\text{ and }(b)$. With \eqref{1au6}, we know
$$f(N(u_5)-\{u_2,u_4\})=\{2_a\}\quad  \mbox{and}\quad  f(N(u_6)-\{u_2\})=\{1_a,1_b\}$$
$$\text{ or }f(N(u_5)-\{u_2,u_4\})=\{2_b\}\quad  \mbox{and}\quad  f(N(u_6)-\{u_2\})=\{1_a,2_a\}.$$  If $f(N(u_5)-\{u_2,u_4\})=\{2_b\}$ and $f(N(u_6)-\{u_2\})=\{1_a,2_a\}$, then we recolor $u_5$ with $2_a$, $u_2$ with $1_a$, and $u$ with $1_b$ to obtain a coloring satisfying $(a)\text{ and }(b)$. Thus, we can assume that 
\begin{equation}\label{u12u13}
f(N(u_5)-\{u_2,u_4\})=\{2_a\}\quad  \mbox{and}\quad  f(N(u_6)-\{u_2\})=\{1_a,1_b\}.
\end{equation}
If $1_b \notin f(N(u_7)-\{u_3\})$, then we can recolor $u_7$ with $1_b$ and it contradicts \eqref{cases}. Thus, we may assume 
\begin{equation}\label{1bu7}
1_b \in f(N(u_7)-\{u_3\}).
\end{equation}
If $1_a \notin f(N(u_8)-\{u_3\})$, then we can recolor $u_8$ with $1_a$ and it contradicts \eqref{u8u9}. If  $1_a \notin f(N(u_{10})-\{u_4\})$, then we can recolor $u_{10}$ with $1_a$ and it contradicts \eqref{cases}. Therefore, we may assume 
\begin{equation}\label{u7u8u10}
1_a \in f(N(u_{10})-\{u_4\})\quad  \mbox{and}\quad   1_a \in f(N(u_8)-\{u_3\}).
\end{equation}
If $2_b \notin f(N(u_7) \cup N(u_8) - \{u_3\})$, then we can recolor $u_3$ with $2_b$ and $u$ with $2_a$ to obtain a coloring satisfying $(a)\text{ and }(b)$.  Thus, we may assume 
\begin{equation}\label{2bu7u8}
2_b \in f(N(u_7) \cup N(u_8) - \{u_3\}).
\end{equation}
By previous arguments, we know that $\{3_a,3_b\} \cap f(\{u_2,u_3,u_4,u_5,u_6,u_7,u_8,u_{10}\}) = \emptyset$.
With \eqref{3a}, \eqref{3b}, and \eqref{u12u13}, we know that $\{3_a,3_b\} \subseteq f(N(u_7) \cup N(u_8) \cup N(u_{10}) - \{u_3,u_4\})$. Moreover, by \eqref{1bu7}, \eqref{u7u8u10}, \eqref{2bu7u8}, and symmetry, we may assume that $$f(N(u_{10})-\{u_4\})=\{1_a,3_b\} \text{ (See Figure~\ref{2-case 2.2.2.})}.$$ But we can recolor $u_{10}$ with $1_b$, $u_4$ with $2_b$, $u_1$ with $1_b$, and $u$ with $1_a$ to obtain a coloring satisfying $(a)\text{ and }(b)$. 


{\bf Case 2.2.3:} $\{u_3u_5, u_4u_6, u_4u_5\} \cap E(G) = \emptyset$ and $u_4u_7 \in E(G)$, i.e., $u_7=u_9$.   If $1_a \notin f(N(u_8)-u_3),$ then we recolor $u_8$ with $1_a$, $u_3$ with $1_b$, and $u$ with $2_a$ to obtain a coloring satisfying $(a) \text{ and } (b).$ Thus, we may assume $1_a \in f(N(u_8)-u_3)$. If $1_a \notin f(N(u_{10})-u_4),$ then we recolor $u_{10}$ with $1_a$, $u_1$ with $2_b$, and $u$ with $1_a$ to obtain a coloring satisfying $(a) \text{ and } (b).$ Thus, we may also assume $1_a \in f(N(u_{10})-u_4)$. If $2_b \notin f(N(u_7) \cup N(u_8)-\{u_3,u_4\})$, then we recolor $u_3$ with $2_b$ and $u$ with $2_a$ to obtain a coloring satisfying $(a)\text{ and }(b)$. With \eqref{3a}, \eqref{3b}, and symmetry, we may assume $f(N(u_7) \cup N(u_8) - \{u_3,u_4\}) = \{1_a,2_b,3_a\}$ and $f(N(u_{10})-u_4) = \{1_a,3_b\}$ (See Figure~\ref{2-case 2.2.3.}). We recolor $u_7$ with $1_b$, $u_4$ with $1_a$, $u_1$ with $1_b$, and $u$ with $1_a$ to obtain a coloring satisfying $(a)\text{ and }(b)$. 

Thus, we may also assume $u_4u_7 \notin E(G)$. 

\medskip
Below we have $\{u_3u_5, u_4u_6, u_4u_5, u_4u_7\} \cap E(G) = \emptyset$. Moreover, by the case (Case 2.2), $$\{u_3u_6,u_4u_8,u_3u_9,u_3u_{10}\} \cap E(G) = \emptyset.$$ Therefore, we also have $|\{u_i: i \in [10]\}|=10$.

\begin{figure}[ht]\label{f4}
\centering
\begin{minipage}[b]{0.4\textwidth}
\begin{tikzpicture}[scale=0.5, transform shape]

\node[circle, draw=white!0,  inner sep=0pt, minimum size=25pt, font=\huge] (11) at (0,1) {$3_a$};
\node[circle, draw=white!0,  inner sep=0pt, minimum size=25pt, font=\huge] (12) at (-4.15,-1) {$1_a$};
\node[circle, draw=white!0,  inner sep=0pt, minimum size=25pt, font=\huge] (14) at (4,-1) {$1_b$};


\node[circle, draw=white!0,  inner sep=0pt, minimum size=25pt, font=\huge] (17) at (-6,-3) {$2_a$};
\node[circle, draw=white!0,  inner sep=0pt, minimum size=25pt, font=\huge] (18) at (-2,-3) {$1_b$};

\node[circle, draw=white!0, inner sep=0pt, minimum size=25pt, font=\huge] (21) at (5.5,-3) {$2_b$};

\node[circle, draw=white!0, inner sep=0pt, minimum size=25pt, font=\huge] (23) at (-7.5,-5) {$1_b$};
\node[circle, draw=white!0, inner sep=0pt, minimum size=25pt, font=\huge] (24) at (2.5,-3) {$1_a$};

\node[circle, draw=white!0, inner sep=0pt, minimum size=25pt, font=\huge] (25) at (-7,-8.5) {$3_a$};

\node[circle, draw=white!0, inner sep=0pt, minimum size=25pt, font=\huge] (27) at (-8,-8.5) {$1_a$};
\node[circle, draw=white!0, inner sep=0pt, minimum size=25pt, font=\huge] (28) at (-0.5,-8.5) {$3_b$};
\node[circle, draw=white!0, inner sep=0pt, minimum size=25pt, font=\huge] (30) at (-4,-5) {$1_a$};

\node[circle, draw=white!0, inner sep=0pt, minimum size=25pt, font=\huge] (32) at (-1.5,-8.5) {$1_a$};
\node[circle, draw=white!0, inner sep=0pt, minimum size=15pt, font=\huge] (33) at (-1,-5) {$2_b$};
\node[circle, draw=white!0, inner sep=0pt, minimum size=15pt, font=\huge] (34) at (-4,-8.5) {$2_b$};

\node[circle, draw=black!80, inner sep=0pt, minimum size=25pt, font=\huge] (1) at (0,0) {$u$};
\node[circle, draw=black!80, inner sep=0pt, minimum size=25pt, font=\huge] (3) at (-4,-2) { $u_1$};
\node[circle, draw=black!80, inner sep=0pt, minimum size=25pt, font=\huge] (4) at (4,-2) {$u_2$};
\node[circle, draw=black!80, inner sep=0pt, minimum size=25pt, font=\huge] (5) at (-6,-4) {$u_3$};
\node[circle, draw=black!80, inner sep=0pt, minimum size=36pt, font=\huge] (6) at (-4,-6) {$u_7$};
\node[circle, draw=black!80, inner sep=0pt, minimum size=36pt, font=\huge] (7) at (-7.5,-6) {$u_8$};
\node[circle, draw=black!80, inner sep=0pt, minimum size=25pt, font=\huge] (8) at (2.5,-4) {$u_5$};
\draw  (1) edge (3);
\draw  (1) edge (4);

\draw  (3) edge (5);
\draw  (4) edge (8);
\draw  (5) edge (6);
\draw  (5) edge (7);

\node (v4) at (-4,-8) {};
\node (v5) at (-8,-8) {};
\node (v6) at (-7,-8) {};

\draw  (6) edge (v4);
\draw  (7) edge (v5);
\draw  (7) edge (v6);
\node[circle, draw=black!80, inner sep=0pt, minimum size=25pt, font=\huge] (29) at (-2,-4) {$u_4$};
\node[circle, draw=black!80, inner sep=0pt, minimum size=25pt, font=\huge] (v7) at (5.5,-4) {$u_6$};
\draw  (3) edge (29);
\draw  (4) edge (v7);

\node[circle, draw=black!80, inner sep=0pt, minimum size=30pt, font=\huge] (v2) at (-1,-6) {$u_{10}$};

\draw  (29) edge (v2);

\node (v10) at (-1.5,-8) {};
\node (v11) at (-0.5,-8) {};

\draw  (v2) edge (v10);
\draw  (v2) edge (v11);

\draw  (6) edge (29);

\end{tikzpicture}

\caption{Case 2.2.3.}
\label{2-case 2.2.3.}
\end{minipage}
\hspace{15mm}
\begin{minipage}[b]{0.4\textwidth}
\begin{tikzpicture}[scale=0.45, transform shape]
\node[circle, draw=white!0,  inner sep=0pt, minimum size=25pt, font=\huge] (11) at (0,1) {$3_a$};
\node[circle, draw=white!0,  inner sep=0pt, minimum size=25pt, font=\huge] (12) at (-4,-1) {$1_a$};
\node[circle, draw=white!0,  inner sep=0pt, minimum size=25pt, font=\huge] (14) at (4,-1) {$1_b$};


\node[circle, draw=white!0,  inner sep=0pt, minimum size=25pt, font=\huge] (17) at (-6,-3) {$2_a$};
\node[circle, draw=white!0,  inner sep=0pt, minimum size=25pt, font=\huge] (18) at (-2,-3) {$1_b$};
\node[circle, draw=white!0,  inner sep=0pt, minimum size=25pt, font=\huge] (19) at (-3.05,-4.85) {$1_a$};

\node[circle, draw=white!0,  inner sep=0pt, minimum size=25pt, font=\huge] (20) at (-1,-4.85) {$2_b$};

\node[circle, draw=white!0, inner sep=0pt, minimum size=25pt, font=\huge] (21) at (5.5,-3) {$2_b$};
\node[circle, draw=white!0, inner sep=0pt, minimum size=25pt, font=\huge] (22) at (-7,-5) {$1_a$};

\node[circle, draw=white!0, inner sep=0pt, minimum size=25pt, font=\huge] (23) at (-5,-4.85) {$1_b$};
\node[circle, draw=white!0, inner sep=0pt, minimum size=25pt, font=\huge] (24) at (2.5,-3) {$1_a$};

\node[circle, draw=white!0, inner sep=0pt, minimum size=25pt, font=\huge] (25) at (-4.5,-7) {$3_a$};
\node[circle, draw=white!0, inner sep=0pt, minimum size=25pt, font=\huge] (26) at (-4.5,-10.5) {$1_b$};

\node[circle, draw=white!0, inner sep=0pt, minimum size=25pt, font=\huge] (28) at (-0.5,-8.5) {$3_b$};
\node[circle, draw=white!0, inner sep=0pt, minimum size=25pt, font=\huge] (30) at (-5.5,-10.5) {$1_a$};
\node[circle, draw=white!0, inner sep=0pt, minimum size=25pt, font=\huge] (31) at (-2.5,-8.5) {$1_b$};
\node[circle, draw=white!0, inner sep=0pt, minimum size=25pt, font=\huge] (32) at (-1.5,-8.5) {$1_a$};
\node[circle, draw=white!0, inner sep=0pt, minimum size=15pt, font=\huge] (33) at (-3.5,-8.5) {$2_b$};
\node[circle, draw=white!0, inner sep=0pt, minimum size=15pt, font=\huge] (34) at (-7.5,-7) {$2_b$};
\node[circle, draw=white!0, inner sep=0pt, minimum size=25pt, font=\huge] (35) at (-3.5,-9.5) {$3_b$};
\node[circle, draw=white!0, inner sep=0pt, minimum size=25pt, font=\huge] (37) at (-0.5,-9.5) {$1_b$};

\node[circle, draw=black!80, inner sep=0pt, minimum size=25pt, font=\huge] (1) at (0,0) {$u$};
\node[circle, draw=black!80, inner sep=0pt, minimum size=25pt, font=\huge] (3) at (-4,-2) { $u_1$};
\node[circle, draw=black!80, inner sep=0pt, minimum size=25pt, font=\huge] (4) at (4,-2) {$u_2$};
\node[circle, draw=black!80, inner sep=0pt, minimum size=25pt, font=\huge] (5) at (-6,-4) {$u_3$};
\node[circle, draw=black!80, inner sep=0pt, minimum size=25pt, font=\huge] (6) at (-7,-6) {$u_7$};
\node[circle, draw=black!80, inner sep=0pt, minimum size=25pt, font=\huge] (7) at (-5,-6) {$u_8$};
\node[circle, draw=black!80, inner sep=0pt, minimum size=25pt, font=\huge] (8) at (2.5,-4) {$u_5$};
\draw  (1) edge (3);
\draw  (1) edge (4);

\draw  (3) edge (5);
\draw  (4) edge (8);
\draw  (5) edge (6);
\draw  (5) edge (7);

\node[circle, draw=black!80, inner sep=0pt, minimum size=25pt, font=\huge] (29) at (-2,-4) {$u_4$};
\node[circle, draw=black!80, inner sep=0pt, minimum size=25pt, font=\huge] (v7) at (5.5,-4) {$u_6$};
\draw  (3) edge (29);
\draw  (4) edge (v7);
\node[circle, draw=black!80, inner sep=0pt, minimum size=25pt, font=\huge] (v1) at (-3,-6) {$u_9$};
\node[circle, draw=black!80, inner sep=0pt, minimum size=25pt, font=\huge] (v2) at (-1,-6) {$u_{10}$};
\draw  (29) edge (v1);
\draw  (29) edge (v2);
\node (v8) at (-3.5,-8) {};
\node (v9) at (-2.5,-8) {};
\node (v10) at (-1.5,-8) {};
\node (v11) at (-0.5,-8) {};
\draw  (v1) edge (v8);
\draw  (v1) edge (v9);
\draw  (v2) edge (v10);
\draw  (v2) edge (v11);

\draw  (6) edge (7);

\node[circle, draw=black!80, inner sep=0pt, minimum size=25pt, font=\huge] (45) at (-5,-8) {$u_{12}$};
\node[circle, draw=black!80, inner sep=0pt, minimum size=25pt, font=\huge] (46) at (-7,-8) {$u_{11}$};

\draw  (6) edge (46);
\draw  (7) edge (45);
\node (v3) at (-5.5,-10) {};
\node (v4) at (-4.5,-10) {};
\draw  (45) edge (v3);
\draw  (45) edge (v4);
\end{tikzpicture}
\caption{Case 2.2.4.}
\label{2-case 2.2.4.}
\end{minipage}
\end{figure}

{\bf Case 2.2.4:} $u_7u_8 \in E(G)$. By \eqref{1a1b'}, $1_b \in f(N(u_9)-\{u_4\})$. If $1_a \notin f(N(u_{10})-\{u_4\})$, then we recolor $u_{10}$ with $1_a$, $u_1$ with $2_b$, and $u$ with $1_a$ to obtain a coloring satisfying $(a)\text{ and }(b)$. Thus, we may assume $1_a \notin f(N(u_{10})-\{u_4\})$. By \eqref{3a} and \eqref{3b}, $\{3_a,3_b\} \subseteq f(N(u_7) \cup N(u_8) \cup N(u_9) \cup N(u_{10}))$. If $\{3_a,3_b\} \subseteq f(N(u_9) \cup N(u_{10}))$, then $f(N(u_9) \cup N(u_{10})-\{u_4\})=\{1_a,1_b,3_a,3_b\}$, $1_b \notin f(N(u_{10})-\{u_4\})$ and $2_b \notin f(N(u_9)-\{u_4\})$. Then, we can recolor $u_{10}$ with $1_b$, $u_4$ with $2_b$, $u_1$ with $1_b$, and $u$ with $1_a$ to obtain a coloring satisfying $(a)\text{ and }(b)$. Thus, by symmetry, we can assume 
\begin{equation}\label{3au7u8}
3_a \in f(N(u_7) \cup N(u_8) - \{u_3\}) \quad \mbox{and} \quad 3_a \notin f(N(u_9) \cup N(u_{10})-u_4).
\end{equation}
If $2_b \notin f(N(u_7) \cup N(u_8) - \{u_3\})$, then we recolor $u_3$ with $2_b$ and $u$ with $2_a$ to obtain a coloring satisfying $(a) \text{ and } (b)$. Thus, we may assume $2_b \in f(N(u_7) \cup N(u_8) - \{u_3\})$. Let $u_{11} \in N(u_7)-\{u_3,u_8\}$ and $u_{12} \in N(u_8)-\{u_3,u_7\}$. We may assume 
\begin{equation}\label{u11u12}
f(u_{11}) = 2_b\quad  \mbox{and}\quad  f(u_{12})=3_a, 
\end{equation}
since, by symmetry, the proof for the case $f(u_{11}) = 3_a$ and $f(u_{12})=2_b$ is similar. 
Note that $3_a \notin f(B(u_1,3)-u_{12})$. If $1_a \notin f(N(u_{12})-\{u_8\})$, then we recolor $u_{12}$ with $1_a$, $u_1$ with $3_a$, and $u$ with $1_a$ to obtain a coloring satisfying $(a)\text{ and }(b)$. If $1_b \notin f(N(u_{12})-\{u_8\})$, then we recolor $u_{12}$ with $1_b$, $u_8$ with $1_a$, $u_7$ with $1_b$, $u_1$ with $3_a$, and $u$ with $1_a$ to obtain a coloring satisfying $(a)\text{ and }(b)$. Thus, we may assume 
\begin{equation}\label{1a1bu12}
f(N(u_{12})-\{u_8\})=\{1_a,1_b\}.
\end{equation}
If $1_b \notin f(N(u_{11})-\{u_7\})$, then we can recolor $u_{11}$ with $1_b$, $u_3$ with $2_b$, and $u$ with $2_a$ to obtain a coloring satisfying $(a)\text{ and }(b)$. Thus, we may assume 
\begin{equation}
1_b \in f(N(u_{11})-\{u_7\}) \text{ (See Figure~\ref{2-case 2.2.4.})}.
\end{equation}
Then, we can recolor $u_8$ with $2_a$, $u_3$ with $1_b$, and $u$ with $2_a$ to obtain a coloring satisfying $(a)\text{ and }(b)$. 

\begin{figure}[ht]\label{f5}
\begin{center}
\begin{tikzpicture}[scale=0.48,transform shape]

\node[circle, draw=white!0,  inner sep=0pt, minimum size=25pt, font=\huge] (11) at (0,1) {$3_a$};
\node[circle, draw=white!0,  inner sep=0pt, minimum size=25pt, font=\huge] (39) at (-4,-1) {$1_a$};
\node[circle, draw=white!0,  inner sep=0pt, minimum size=25pt, font=\huge] (14) at (4,-1) {$1_b$};


\node[circle, draw=white!0,  inner sep=0pt, minimum size=25pt, font=\huge] (17) at (-6.5,-3) {$2_a$};
\node[circle, draw=white!0,  inner sep=0pt, minimum size=25pt, font=\huge] (30) at (-1.5,-3) {$1_b$};
\node[circle, draw=white!0,  inner sep=0pt, minimum size=25pt, font=\huge] (19) at (-2.5,-5) {$1_a$};

\node[circle, draw=white!0,  inner sep=0pt, minimum size=25pt, font=\huge] (20) at (-0.5,-5) {$2_b$};

\node[circle, draw=white!0, inner sep=0pt, minimum size=25pt, font=\huge] (21) at (5.5,-3) {$2_b$};
\node[circle, draw=white!0, inner sep=0pt, minimum size=25pt, font=\huge] (22) at (-7.5,-5) {$1_a$};

\node[circle, draw=white!0, inner sep=0pt, minimum size=25pt, font=\huge] (23) at (-5.5,-5) {$1_b$};
\node[circle, draw=white!0, inner sep=0pt, minimum size=25pt, font=\huge] (24) at (2.5,-3) {$1_a$};

\node[circle, draw=white!0, inner sep=0pt, minimum size=25pt, font=\huge] (25) at (-8,-8.5) {$3_a$};

\node[circle, draw=white!0, inner sep=0pt, minimum size=15pt, font=\huge] (34) at (-5.5,-8.5) {$2_b$};
\node[circle, draw=white!0, inner sep=0pt, minimum size=25pt, font=\huge] (35) at (-2.5,-8.5) {$3_b$};
\node[circle, draw=white!0, inner sep=0pt, minimum size=25pt, font=\huge] (37) at (-7,-8.5) {$1_b$};

\node[circle, draw=black!80, inner sep=0pt, minimum size=25pt, font=\huge] (1) at (0,0) {$u$};
\node[circle, draw=black!80, inner sep=0pt, minimum size=25pt, font=\huge] (3) at (-4,-2) { $u_1$};
\node[circle, draw=black!80, inner sep=0pt, minimum size=25pt, font=\huge] (4) at (4,-2) {$u_2$};
\node[circle, draw=black!80, inner sep=0pt, minimum size=25pt, font=\huge] (5) at (-6.5,-4) {$u_3$};
\node[circle, draw=black!80, inner sep=0pt, minimum size=25pt, font=\huge] (6) at (-7.5,-6) {$u_7$};
\node[circle, draw=black!80, inner sep=0pt, minimum size=25pt, font=\huge] (7) at (-5.5,-6) {$u_8$};
\node[circle, draw=black!80, inner sep=0pt, minimum size=25pt, font=\huge] (8) at (2.5,-4) {$u_5$};
\draw  (1) edge (3);
\draw  (1) edge (4);

\draw  (3) edge (5);
\draw  (4) edge (8);
\draw  (5) edge (6);
\draw  (5) edge (7);

\node[circle, draw=black!80, inner sep=0pt, minimum size=25pt, font=\huge] (29) at (-1.5,-4) {$u_4$};
\node[circle, draw=black!80, inner sep=0pt, minimum size=25pt, font=\huge] (v7) at (5.5,-4) {$u_6$};
\draw  (3) edge (29);
\draw  (4) edge (v7);
\node[circle, draw=black!80, inner sep=0pt, minimum size=25pt, font=\huge] (v1) at (-2.5,-6) {$u_9$};
\node[circle, draw=black!80, inner sep=0pt, minimum size=25pt, font=\huge] (v2) at (-0.5,-6) {$u_{10}$};
\draw  (29) edge (v1);
\draw  (29) edge (v2);
\node (v8) at (-2.5,-8) {};

\node (v10) at (-1,-8) {};
\node (v11) at (0,-8) {};
\draw  (v1) edge (v8);

\draw  (v2) edge (v10);
\draw  (v2) edge (v11);

\node(45) at (-5.5,-8) {};
\node(46) at (-8,-8) {};

\draw  (6) edge (46);
\draw  (7) edge (45);

\draw  (7) edge (v1);

\node (v3) at (-7,-8) {};
\draw  (6) edge (v3);

\end{tikzpicture}
\caption{Case 2.2.5.}
\label{2-case 2.2.5.}
\end{center}
\end{figure}
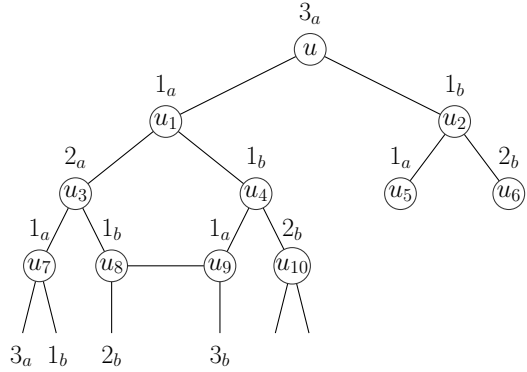

{\bf Case 2.2.5:} $u_7u_8 \notin E(G), u_8u_9 \in E(G)$.  Similarly to \eqref{1au10} and \eqref{1bu7}, we may assume 
\begin{equation}\label{u7u10}
1_a \in f(N(u_{10})-\{u_4\})\quad  \mbox{and}\quad  1_b \in f(N(u_7)-\{u_3\}).
\end{equation}
If $2_b \notin f(N(u_7) \cup N(u_8) - \{u_3\})$, then we recolor $u_3$ with $2_b$ and $u$ with $2_a$ to obtain a coloring satisfying $(a)\text{ and }(b)$. Thus, we may assume $2_b \in f(N(u_7) \cup N(u_8) - \{u_3\})$. If $1_b \notin f(N(u_{10})-\{u_4\})$ and $2_b \notin f(N(u_9)-\{u_4\})$, then we can recolor $u_{10}$ with $1_b$, $u_4$ with $2_b$, $u_1$ with $1_b$, and $u$ with $1_a$ to obtain a coloring satisfying $(a)\text{ and }(b)$. From \eqref{3a} and \eqref{3b}, we know that $f(N(u_8) \cup N(u_9)-\{u_3,u_4\}) \subseteq \{2_b,3_a,3_b\}$ (See Figure~\ref{2-case 2.2.5.}). But it contradicts \eqref{1a1b'}. 

Therefore, we may assume $u_8u_9 \notin E(G)$.

{ \bf Case 2.2.6:} $u_7u_8 \notin E(G), u_8u_9 \notin E(G)$.
If $|N(u_7)|=|N(u_8)|=|N(u_9)|=|N(u_{10})|=3,$ then we let $$\{u_{11},u_{12}\} \subseteq N(u_7)-\{u_3\}, \quad \{u_{13},u_{14}\} \subseteq N(u_8)-\{u_3\}, \quad \{u_{15},u_{16}\} \subseteq N(u_9)-\{u_4\},$$$$\quad  \mbox{and}\quad  \{u_{17},u_{18}\} \subseteq N(u_{10})-\{u_4\}.$$ 
It is possible that $|\{u_i:i \in[18]-[10]\}| \neq 8$ or $\{u_5,u_6\} \cap \{u_i: i \in [18]-[10]\} \neq \emptyset$, but this will not affect the proof below.

Similarly to \eqref{1bu9}, \eqref{1bu5}, \eqref{1au10}, \eqref{u8}, we may assume
\begin{equation}\label{1a1b226}
f(u_{12})=f(u_{16})=1_b \quad  \mbox{and}\quad   f(u_{13})=f(u_{17})=1_a.
\end{equation}
Similarly to \eqref{2bu7u8} and \eqref{3au7u8}, we may assume 
\begin{equation}\label{2b3au7u8}
\{2_b,3_a\} \subseteq f(N(u_7) \cup N(u_8) - \{u_3\}).
\end{equation}
If $1_b \notin f(N(u_{10})-\{u_4\})$ and $2_b \notin f(N(u_9)-\{u_4\})$, then we can recolor $u_{10}$ with $1_b$, $u_4$ with $2_b$, $u_1$ with $1_b$, and $u$ with $1_a$ to obtain a coloring satisfying $(a)\text{ and }(b)$. With \eqref{3b}, we may assume

\begin{equation}\label{u9u10}
\mbox{\em either
$f(u_{15})=3_b$ and  $f(u_{18})=1_b$ or $f(u_{15})=2_b$ and  $f(u_{18})=3_b$.}
\end{equation}

If $|N(u_{11})|=|N(u_{12})|=|N(u_{13})|=|N(u_{14})|=3$, then we let $\{u_{19},u_{20}\} \subseteq N(u_{11})$, $\{u_{21},u_{22}\} \subseteq N(u_{12})$, $\{u_{23},u_{24}\} \subseteq N(u_{13})$, $\{u_{25},u_{26}\} \subseteq N(u_{14})$. 

By \eqref{2b3au7u8}, we have 

\begin{equation}\label{cases26}
\mbox{\em either 
$f(u_{11})=2_b$ and $f(u_{14})=3_a$ or $f(u_{11})=3_a$  and $f(u_{14})=2_b$.
}
\end{equation}


\begin{figure}[ht]\label{f6}
\begin{center}
\begin{tikzpicture}[scale=0.35,transform shape]

\node[circle, draw=white!0,   inner sep=0pt, minimum size=25pt, font=\huge] (51) at (-5.8,-12.35) {$2_b$};
\node[circle, draw=white!0,   inner sep=0pt, minimum size=25pt, font=\huge] (52) at (-12.2,-12.35) {$2_a$};

\node[circle, draw=white!0,   inner sep=0pt, minimum size=25pt, font=\huge] (54) at (-8.5,-17.5) {$2_a$};

\node[circle, draw=white!0,   inner sep=0pt, minimum size=15pt, font=\huge] (58) at (3,-4.85) {$2_b$};
\node[circle, draw=white!0,   inner sep=0pt, minimum size=25pt, font=\huge] (59) at (-11.5,-17.5) {$1_b$};
\node[circle, draw=white!0,   inner sep=0pt, minimum size=25pt, font=\huge] (60) at (-7.5,-17.5) {$1_b$};

\node[circle, draw=white!0,   inner sep=0pt, minimum size=25pt, font=\huge] (11) at (1,1) {$3_a$};
\node[circle, draw=white!0,   inner sep=0pt, minimum size=25pt, font=\huge] (12) at (-12.5,-17.5) {$3_a$};
\node[circle, draw=white!0,   inner sep=0pt, minimum size=25pt, font=\huge] (14) at (-18,-14.65) {$1_b$};
\node[circle, draw=white!0,   inner sep=0pt, minimum size=25pt, font=\huge] (36) at (-10.5,-17.5) {$1_a$};

\node[circle, draw=white!0,   inner sep=0pt, minimum size=25pt, font=\huge] (18) at (-0.5,-3) {$1_b$};
\node[circle, draw=white!0,   inner sep=0pt, minimum size=25pt, font=\huge] (19) at (-4,-4.85) {$1_a$};

\node[circle, draw=white!0,   inner sep=0pt, minimum size=25pt, font=\huge] (20) at (-13,-3) {$2_a$};

\node[circle, draw=white!0,   inner sep=0pt, minimum size=25pt, font=\huge] (21) at (11,-3) {$2_b$};
\node[circle, draw=white!0,   inner sep=0pt, minimum size=25pt, font=\huge] (22) at (-19,-8.35) {$2_b$};

\node[circle, draw=white!0,   inner sep=0pt, minimum size=25pt, font=\huge] (23) at (-5.5,-17.5) {$1_b$};
\node[circle, draw=white!0,   inner sep=0pt, minimum size=25pt, font=\huge] (24) at (-6.5,-17.5) {$1_a$};

\node[circle, draw=white!0,   inner sep=0pt, minimum size=25pt, font=\huge] (25) at (-7,-8.35) {$3_a$};
\node[circle, draw=white!0,   inner sep=0pt, minimum size=25pt, font=\huge] (26) at (-9.8,-12.35) {$1_b$};
\node[circle, draw=white!0,   inner sep=0pt, minimum size=25pt, font=\huge] (27) at (-8.2,-12.35) {$1_a$};
\node[circle, draw=white!0,   inner sep=0pt, minimum size=25pt, font=\huge] (28) at (-9.5,-17.5) {$3_b$};
\node[circle, draw=white!0,   inner sep=0pt, minimum size=25pt, font=\huge] (30) at (-14.5,-14.65) {$3_b$};
\node[circle, draw=white!0,   inner sep=0pt, minimum size=25pt, font=\huge] (61) at (-14.5,-15.5) {$2_b$};
\node[circle, draw=white!0,   inner sep=0pt, minimum size=25pt, font=\huge] (31) at (4,-10.65) {$1_b$};
\node[circle, draw=white!0,   inner sep=0pt, minimum size=25pt, font=\huge] (64) at (4,-11.5) {$3_b$};
\node[circle, draw=white!0,   inner sep=0pt, minimum size=25pt, font=\huge] (32) at (-16.5,-14.65) {$1_a$};

\node[circle, draw=white!0,   inner sep=0pt, minimum size=25pt, font=\huge] (35) at (-5,-10.65) {$3_b$};
\node[circle, draw=white!0,   inner sep=0pt, minimum size=25pt, font=\huge] (63) at (-5,-11.5) {$2_b$};
\node[circle, draw=white!0,   inner sep=0pt, minimum size=25pt, font=\huge] (37) at (-15.5,-8.35) {$1_b$};
\node[circle, draw=white!0,   inner sep=0pt, minimum size=25pt, font=\huge] (38) at (-20,-14.65) {$1_a$};
\node[circle, draw=white!0,   inner sep=0pt, minimum size=25pt, font=\huge] (60) at (-20,-15.5) {$3_b$};
\node[circle, draw=white!0,   inner sep=0pt, minimum size=25pt, font=\huge] (39) at (-9,-4.85) {$1_b$};
\node[circle, draw=white!0,   inner sep=0pt, minimum size=25pt, font=\huge] (40) at (6,-6.5) {$2_a$};
\node[circle, draw=white!0,   inner sep=0pt, minimum size=25pt, font=\huge] (41) at (2,-10.65) {$1_a$};
\node[circle, draw=white!0,   inner sep=0pt, minimum size=25pt, font=\huge] (42) at (-11,-8.35) {$1_a$};
\node[circle, draw=white!0,   inner sep=0pt, minimum size=25pt, font=\huge] (43) at (-17,-4.85) {$1_a$};
\node[circle, draw=white!0,   inner sep=0pt, minimum size=25pt, font=\huge] (44) at (10,-6.5) {$1_a$};
\node[circle, draw=white!0,   inner sep=0pt, minimum size=25pt, font=\huge] (45) at (7,-3) {$1_a$};
\node[circle, draw=white!0,   inner sep=0pt, minimum size=25pt, font=\huge] (46) at (-7,-0.9) {$1_a$};
\node[circle, draw=white!0,   inner sep=0pt, minimum size=25pt, font=\huge] (47) at (-3,-10.65) {$1_b$};
\node[circle, draw=white!0,   inner sep=0pt, minimum size=25pt, font=\huge] (48) at (12,-6.5) {$1_b$};
\node[circle, draw=white!0,   inner sep=0pt, minimum size=25pt, font=\huge] (49) at (8,-6.5) {$1_b$};
\node[circle, draw=white!0,   inner sep=0pt, minimum size=25pt, font=\huge] (50) at (9,-0.9) {$1_b$};

\node[circle, draw=black!80,   inner sep=0pt, minimum size=25pt, font=\huge] (1) at (1,0) {$u$};
\node[circle, draw=black!80,   inner sep=0pt, minimum size=25pt, font=\huge] (3) at (-7,-2) { $u_1$};
\node[circle, draw=black!80,   inner sep=0pt, minimum size=25pt, font=\huge] (4) at (9,-2) {$u_2$};
\node[circle, draw=black!80,   inner sep=0pt, minimum size=25pt, font=\huge] (5) at (-13,-4) {$u_3$};
\node[circle, draw=black!80,   inner sep=0pt, minimum size=36pt, font=\huge] (6) at (-17,-6) {$u_7$};
\node[circle, draw=black!80,   inner sep=0pt, minimum size=36pt, font=\huge] (7) at (-9,-6) {$u_8$};
\node[circle, draw=black!80,   inner sep=0pt, minimum size=25pt, font=\huge] (8) at (7,-4) {$u_5$};
\draw  (1) edge (3);
\draw  (1) edge (4);

\draw  (3) edge (5);
\draw  (4) edge (8);
\draw  (5) edge (6);
\draw  (5) edge (7);

\node[circle, draw=black!80,   inner sep=0pt, minimum size=30pt, font=\huge] (v3) at (-19,-9.5) {$u_{11}$};
\node[circle, draw=black!80,   inner sep=0pt, minimum size=30pt, font=\huge] (v4) at (-15.5,-9.5) {$u_{12}$};
\node[circle, draw=black!80,   inner sep=0pt, minimum size=30pt, font=\huge] (v5) at (-11,-9.5) {$u_{13}$};
\node[circle, draw=black!80,   inner sep=0pt, minimum size=30pt, font=\huge] (v6) at (-7,-9.5) {$u_{14}$};
\draw  (6) edge (v3);
\draw  (6) edge (v4);
\draw  (7) edge (v5);
\draw  (7) edge (v6);
\node[circle, draw=black!80,   inner sep=0pt, minimum size=25pt, font=\huge] (29) at (-0.5,-4) {$u_4$};
\node[circle, draw=black!80,   inner sep=0pt, minimum size=25pt, font=\huge] (v7) at (11,-4) {$u_6$};
\draw  (3) edge (29);
\draw  (4) edge (v7);
\node[circle, draw=black!80,   inner sep=0pt, minimum size=36pt, font=\huge] (v1) at (-4,-6) {$u_9$};
\node[circle, draw=black!80,   inner sep=0pt, minimum size=36pt, font=\huge] (v2) at (3,-6) {$u_{10}$};
\draw  (29) edge (v1);
\draw  (29) edge (v2);
\node[circle, draw=black!80,   inner sep=0pt, minimum size=30pt, font=\huge] (v8) at (-5,-9.5) {$u_{15}$};
\node[circle, draw=black!80,   inner sep=0pt, minimum size=30pt, font=\huge] (v9) at (-3,-9.5) {$u_{16}$};
\node[circle, draw=black!80,   inner sep=0pt, minimum size=30pt, font=\huge] (v10) at (2,-9.5) {$u_{17}$};
\node[circle, draw=black!80,   inner sep=0pt, minimum size=30pt, font=\huge] (v11) at (4,-9.5) {$u_{18}$};
\draw  (v1) edge (v8);
\draw  (v1) edge (v9);
\draw  (v2) edge (v10);
\draw  (v2) edge (v11);
\node[circle, draw=black!80,   inner sep=0pt, minimum size=30pt, font=\huge] (v12) at (-20,-13.5) {$u_{19}$};
\node[circle, draw=black!80,   inner sep=0pt, minimum size=30pt, font=\huge] (v14) at (-18,-13.5) {$u_{20}$};
\node[circle, draw=black!80,   inner sep=0pt, minimum size=30pt, font=\huge] (v13) at (-16.5,-13.5) {$u_{21}$};
\node[circle, draw=black!80,   inner sep=0pt, minimum size=30pt, font=\huge] (v15) at (-14.5,-13.5) {$u_{22}$};
\node[circle, draw=black!80,   inner sep=0pt, minimum size=30pt, font=\huge] (v16) at (-12,-13.5) {$u_{23}$};
\node[circle, draw=black!80,   inner sep=0pt, minimum size=30pt, font=\huge] (v17) at (-10,-13.5) {$u_{24}$};
\node[circle, draw=black!80,   inner sep=0pt, minimum size=30pt, font=\huge] (v18) at (-8,-13.5) {$u_{25}$};
\node[circle, draw=black!80,   inner sep=0pt, minimum size=30pt, font=\huge] (v19) at (-6,-13.5) {$u_{26}$};
\node (v20) at (-12.5,-17) {};
\node (v22) at (-10.5,-17) {};
\node (v21) at (-11.5,-17) {};
\node (v23) at (-9.5,-17) {};
\node (v24) at (-8.5,-17) {};
\node (v25) at (-7.5,-17) {};
\node (v26) at (-6.5,-17) {};
\node (v27) at (-5.5,-17) {};
\draw  (v3) edge (v12);
\draw  (v4) edge (v13);
\draw  (v3) edge (v14);
\draw  (v4) edge (v15);
\draw  (v5) edge (v16);
\draw  (v5) edge (v17);
\draw  (v6) edge (v18);
\draw  (v6) edge (v19);
\draw  (v16) edge (v20);
\draw  (v16) edge (v21);
\draw  (v17) edge (v22);
\draw  (v17) edge (v23);
\draw  (v18) edge (v24);
\draw  (v18) edge (v25);
\draw  (v19) edge (v26);
\draw  (v19) edge (v27);

\node (v28) at (6,-6) {};
\node (v29) at (8,-6) {};
\node (v30) at (10,-6) {};
\node (v31) at (12,-6) {};
\draw  (8) edge (v28);
\draw  (8) edge (v29);
\draw  (v7) edge (v30);
\draw  (v7) edge (v31);

\end{tikzpicture}
\caption{Case 2.2.6.1.}
\label{2-case 2.2.6.1.}
\end{center}
\end{figure}
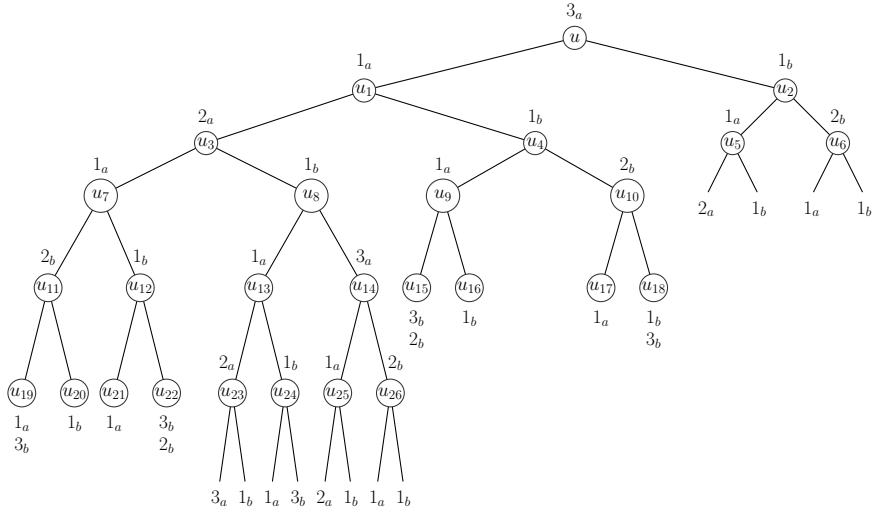


{\bf Case 2.2.6.1:} $f(u_{11})=2_b$ and $f(u_{14}) = 3_a$. If $1_b \notin f(N(u_{13})-\{u_8\})$, then we can recolor $u_{13}$ with $1_b$, $u_8$ with $1_a$, $u_3$ with $1_b$, and $u$ with $2_a$ to obtain a coloring satisfying $(a)\text{ and }(b)$. If $2_b \notin f(N(u_{13}) \cup N(u_{14})-\{u_8\})$, then we can recolor $u_8$ with $2_b$, $u_3$ with $1_b$, and $u$ with $2_a$ to obtain a coloring satisfying $(a)\text{ and }(b)$. Thus, we may assume 
\begin{equation}\label{2bu13u14}
2_b \in f(N(u_{13}) \cup N(u_{14})-\{u_8\}). 
\end{equation}
If $2_a \notin f(N(u_{13}) \cup N(u_{14})-\{u_8\})$, then we can recolor $u_8$ with $2_a$, $u_3$ with $1_b$, and $u$ with $2_a$ to obtain a coloring satisfying $(a)\text{ and }(b)$. Thus, we may also assume
\begin{equation}\label{2au13u14}
2_a \in f(N(u_{13}) \cup N(u_{14})-\{u_8\}). 
\end{equation}
If $1_b \notin f(N(u_{11})-\{u_7\})$, then we can recolor $u_{11}$ with $1_b$ and it contradicts \eqref{cases26}. Similarly, $1_a \in f(N(u_{14})-\{u_8\})$. If $1_a \notin f(N(u_{12})-\{u_7\})$, then we can recolor $u_{12}$ with $1_a$, $u_7$ with $1_b$, and it contradicts \eqref{cases}. Similarly, $1_b \in f(N(u_{13})-\{u_8\})$. Thus, we may assume 
\begin{equation}\label{u24u25}
|N(u_{13})|=|N(u_{14})|=3, f(u_{20})=f(u_{24})=1_b,\quad  \mbox{and}\quad  f(u_{21})=f(u_{25})=1_a.
\end{equation}
Furthermore, by \eqref{2bu13u14} and \eqref{2au13u14}, we assume
\begin{equation}\label{u23u26}
f(u_{23})=2_a\quad  \mbox{and}\quad  f(u_{26})=2_b,
\end{equation} 
since the argument for $f(u_{23})=2_b$ and $f(u_{26})=2_a$ is similar. 
If $\{1_a,1_b\} \neq f(N(u_{26})-\{u_{14}\})$, then we can recolor $u_{26}$ with a color $x \in f(N(u_{26})-\{u_{14}\})-\{1_a,1_b\}$, $u_8$ with $2_b$, $u_3$ with $1_b$, and $u$ with $2_a$ to obtain a coloring satisfying $(a)\text{ and }(b)$. Thus, we may assume 
\begin{equation}\label{u26}
f(N(u_{26})-\{u_{14}\})=\{1_a,1_b\}.
\end{equation}
If $1_b \notin f(N(u_{25})-\{u_{14}\})$, then we can recolor $u_{25}$ with $1_b$, $u_{14}$ with $1_a$, and it contradicts \eqref{cases26}. Thus, we may assume 
\begin{equation}
1_b \in f(N(u_{25})-\{u_{14}\}).
\end{equation}
If $f(u_{19}) \neq 1_a$ and $f(u_{22}) \neq 2_b$, then we can recolor $u_{11}$ with $1_a$, $u_7$ with $2_b$, $u_3$ with $1_a$, $u_1$ with $2_a$, and $u$ with $1_a$ to obtain a coloring satisfying $(a)\text{ and }(b)$. If $3_b \notin f(N(u_{11}) \cup N(u_{12})-\{u_7\})$, then we can recolor $u_3$ with $3_b$ and $u$ with $2_a$ to obtain a coloring satisfying $(a)\text{ and }(b)$. Thus, we can assume
\begin{equation}\label{u19u22}
\mbox{\em either
$f(u_{19})=1_a$ and  $f(u_{22})=3_b$ or $f(u_{19})=3_b$ and $f(u_{22})=2_b$.
}
\end{equation}

If $2_a \notin f(N(u_{25}) \cup N(u_{26}) - \{u_{14}\})$, then by \eqref{u19u22}, we can recolor $u_{14}$ with $2_a$, $u_3$ with $3_a$, and $u$ with $2_a$ to obtain a coloring satisfying $(a)\text{ and }(b)$. With \eqref{u26}, we may assume 
\begin{equation}
2_a \in f(N(u_{25})-\{u_{14}\}).
\end{equation}
Similarly to \eqref{u24u25}, we may assume 
\begin{equation}\label{u23u24}
1_a \in f(N(u_{24})-\{u_{13}\})\quad  \mbox{and}\quad  1_b \in f(N(u_{23})-\{u_{13}\}).
\end{equation}
If $\{3_a,3_b\} \nsubseteq f(N(u_{23}) \cup N(u_{24}))$, then we can recolor $u_8$ with a color $x \in f(N(u_{23}) \cup N(u_{24})) - \{3_a,3_b\}$, $u_{14}$ with $1_b$, $u_3$ with $1_b$, and $u$ with $2_a$ to obtain a coloring satisfying $(a)\text{ and }(b)$. Therefore, 
$$ f(N(u_{23}) \cup N(u_{24})-\{u_{13}\})=\{1_a,1_b,3_a,3_b\}\quad  \mbox{and}\quad  2_b \notin f(B(u_{13})) \text{ (See Figure~\ref{2-case 2.2.6.1.})}.$$
We recolor $u_{13}$ with $2_b$, $u_8$ with $1_a$, $u_3$ with $1_b$, and $u$ with $2_a$ to obtain a coloring satisfying $(a)\text{ and }(b)$.

\begin{figure}[ht]\label{f7}
\begin{center}
\begin{tikzpicture}[scale=0.35,transform shape]

\node[circle, draw=white!0,   inner sep=0pt, minimum size=25pt, font=\huge] (11) at (1,1) {$3_a$};

\node[circle, draw=white!0,   inner sep=0pt, minimum size=25pt, font=\huge] (14) at (-18,-14.65) {$1_b$};


\node[circle, draw=white!0,   inner sep=0pt, minimum size=25pt, font=\huge] (18) at (-0.5,-3) {$1_b$};
\node[circle, draw=white!0,   inner sep=0pt, minimum size=25pt, font=\huge] (19) at (-4,-4.85) {$1_a$};

\node[circle, draw=white!0,   inner sep=0pt, minimum size=25pt, font=\huge] (20) at (-13,-3) {$2_a$};

\node[circle, draw=white!0,   inner sep=0pt, minimum size=25pt, font=\huge] (21) at (11,-3) {$2_b$};
\node[circle, draw=white!0,   inner sep=0pt, minimum size=25pt, font=\huge] (22) at (-19,-8.35) {$3_a$};
\node[circle, draw=white!0,   inner sep=0pt, minimum size=25pt, font=\huge] (25) at (-7,-8.35) {$2_b$};

\node[circle, draw=white!0,   inner sep=0pt, minimum size=25pt, font=\huge] (30) at (-14.5,-14.65) {$3_b$};

\node[circle, draw=white!0,   inner sep=0pt, minimum size=25pt, font=\huge] (71) at (-14.5,-15.5) {$2_b$};

\node[circle, draw=white!0,   inner sep=0pt, minimum size=25pt, font=\huge] (31) at (4,-10.65) {$1_b$};
\node[circle, draw=white!0,   inner sep=0pt, minimum size=25pt, font=\huge] (77) at (4,-11.5) {$3_b$};

\node[circle, draw=white!0,   inner sep=0pt, minimum size=25pt, font=\huge] (32) at (-16.5,-14.65) {$1_a$};

\node[circle, draw=white!0,   inner sep=0pt, minimum size=25pt, font=\huge] (35) at (-5,-10.65) {$3_b$};
\node[circle, draw=white!0,   inner sep=0pt, minimum size=25pt, font=\huge] (75) at (-5,-11.5) {$2_b$};
\node[circle, draw=white!0,   inner sep=0pt, minimum size=25pt, font=\huge] (37) at (-15.5,-8.35) {$1_b$};
\node[circle, draw=white!0,   inner sep=0pt, minimum size=25pt, font=\huge] (38) at (-20,-14.65) {$2_b$};
\node[circle, draw=white!0,   inner sep=0pt, minimum size=25pt, font=\huge] (70) at (-20,-15.5) {$3_b$};
\node[circle, draw=white!0,   inner sep=0pt, minimum size=25pt, font=\huge] (39) at (-9,-4.85) {$1_b$};
\node[circle, draw=white!0,   inner sep=0pt, minimum size=25pt, font=\huge] (40) at (6,-6.5) {$2_a$};
\node[circle, draw=white!0,   inner sep=0pt, minimum size=25pt, font=\huge] (41) at (2,-10.65) {$1_a$};
\node[circle, draw=white!0,   inner sep=0pt, minimum size=25pt, font=\huge] (42) at (-11,-8.35) {$1_a$};
\node[circle, draw=white!0,   inner sep=0pt, minimum size=25pt, font=\huge] (43) at (-17,-4.85) {$1_a$};
\node[circle, draw=white!0,   inner sep=0pt, minimum size=25pt, font=\huge] (44) at (10,-6.5) {$1_a$};
\node[circle, draw=white!0,   inner sep=0pt, minimum size=25pt, font=\huge] (45) at (7,-3) {$1_a$};
\node[circle, draw=white!0,   inner sep=0pt, minimum size=25pt, font=\huge] (46) at (-7,-1) {$1_a$};
\node[circle, draw=white!0,   inner sep=0pt, minimum size=25pt, font=\huge] (47) at (-3,-10.65) {$1_b$};
\node[circle, draw=white!0,   inner sep=0pt, minimum size=25pt, font=\huge] (48) at (12,-6.5) {$1_b$};
\node[circle, draw=white!0,   inner sep=0pt, minimum size=25pt, font=\huge] (49) at (8,-6.5) {$1_b$};
\node[circle, draw=white!0,   inner sep=0pt, minimum size=25pt, font=\huge] (50) at (9,-1) {$1_b$};

\node[circle, draw=white!0,   inner sep=0pt, minimum size=15pt, font=\huge] (58) at (3,-4.85) {$2_b$};
\node[circle, draw=white!0,   inner sep=0pt, minimum size=15pt, font=\huge] (59) at (-6,-14.65) {$1_b$};
\node[circle, draw=white!0,   inner sep=0pt, minimum size=15pt, font=\huge] (72) at (-6,-15.5) {$2_a$};
\node[circle, draw=white!0,   inner sep=0pt, minimum size=15pt, font=\huge] (60) at (-8,-14.65) {$1_a$};
\node[circle, draw=white!0,   inner sep=0pt, minimum size=15pt, font=\huge] (61) at (-10,-14.65) {$1_b$};
\node[circle, draw=white!0,   inner sep=0pt, minimum size=15pt, font=\huge] (62) at (-12,-14.65) {$2_a$};
\node[circle, draw=white!0,   inner sep=0pt, minimum size=15pt, font=\huge] (68) at (-12,-15.5) {$2_b$};

\node[circle, draw=black!80,   inner sep=0pt, minimum size=25pt, font=\huge] (1) at (1,0) {$u$};
\node[circle, draw=black!80,   inner sep=0pt, minimum size=25pt, font=\huge] (3) at (-7,-2) { $u_1$};
\node[circle, draw=black!80,   inner sep=0pt, minimum size=25pt, font=\huge] (4) at (9,-2) {$u_2$};
\node[circle, draw=black!80,   inner sep=0pt, minimum size=25pt, font=\huge] (5) at (-13,-4) {$u_3$};
\node[circle, draw=black!80,   inner sep=0pt, minimum size=36pt, font=\huge] (6) at (-17,-6) {$u_7$};
\node[circle, draw=black!80,   inner sep=0pt, minimum size=36pt, font=\huge] (7) at (-9,-6) {$u_8$};
\node[circle, draw=black!80,   inner sep=0pt, minimum size=25pt, font=\huge] (8) at (7,-4) {$u_5$};
\draw  (1) edge (3);
\draw  (1) edge (4);

\draw  (3) edge (5);
\draw  (4) edge (8);
\draw  (5) edge (6);
\draw  (5) edge (7);

\node[circle, draw=black!80,   inner sep=0pt, minimum size=30pt, font=\huge] (v3) at (-19,-9.5) {$u_{11}$};
\node[circle, draw=black!80,   inner sep=0pt, minimum size=30pt, font=\huge] (v4) at (-15.5,-9.5) {$u_{12}$};
\node[circle, draw=black!80,   inner sep=0pt, minimum size=30pt, font=\huge] (v5) at (-11,-9.5) {$u_{13}$};
\node[circle, draw=black!80,   inner sep=0pt, minimum size=30pt, font=\huge] (v6) at (-7,-9.5) {$u_{14}$};
\draw  (6) edge (v3);
\draw  (6) edge (v4);
\draw  (7) edge (v5);
\draw  (7) edge (v6);
\node[circle, draw=black!80,   inner sep=0pt, minimum size=25pt, font=\huge] (29) at (-0.5,-4) {$u_4$};
\node[circle, draw=black!80,   inner sep=0pt, minimum size=25pt, font=\huge] (v7) at (11,-4) {$u_6$};
\draw  (3) edge (29);
\draw  (4) edge (v7);
\node[circle, draw=black!80,   inner sep=0pt, minimum size=36pt, font=\huge] (v1) at (-4,-6) {$u_9$};
\node[circle, draw=black!80,   inner sep=0pt, minimum size=36pt, font=\huge] (v2) at (3,-6) {$u_{10}$};
\draw  (29) edge (v1);
\draw  (29) edge (v2);
\node[circle, draw=black!80,   inner sep=0pt, minimum size=30pt, font=\huge] (v8) at (-5,-9.5) {$u_{15}$};
\node[circle, draw=black!80,   inner sep=0pt, minimum size=30pt, font=\huge] (v9) at (-3,-9.5) {$u_{16}$};
\node[circle, draw=black!80,   inner sep=0pt, minimum size=30pt, font=\huge] (v10) at (2,-9.5) {$u_{17}$};
\node[circle, draw=black!80,   inner sep=0pt, minimum size=30pt, font=\huge] (v11) at (4,-9.5) {$u_{18}$};
\draw  (v1) edge (v8);
\draw  (v1) edge (v9);
\draw  (v2) edge (v10);
\draw  (v2) edge (v11);
\node[circle, draw=black!80,   inner sep=0pt, minimum size=30pt, font=\huge] (v12) at (-20,-13.5) {$u_{19}$};
\node[circle, draw=black!80,   inner sep=0pt, minimum size=30pt, font=\huge] (v14) at (-18,-13.5) {$u_{20}$};
\node[circle, draw=black!80,   inner sep=0pt, minimum size=30pt, font=\huge] (v13) at (-16.5,-13.5) {$u_{21}$};
\node[circle, draw=black!80,   inner sep=0pt, minimum size=30pt, font=\huge] (v15) at (-14.5,-13.5) {$u_{22}$};
\node[circle, draw=black!80,   inner sep=0pt, minimum size=30pt, font=\huge] (v16) at (-12,-13.5) {$u_{23}$};
\node[circle, draw=black!80,   inner sep=0pt, minimum size=30pt, font=\huge] (v17) at (-10,-13.5) {$u_{24}$};
\node[circle, draw=black!80,   inner sep=0pt, minimum size=30pt, font=\huge] (v18) at (-8,-13.5) {$u_{25}$};
\node[circle, draw=black!80,   inner sep=0pt, minimum size=30pt, font=\huge] (v19) at (-6,-13.5) {$u_{26}$};

\draw  (v3) edge (v12);
\draw  (v4) edge (v13);
\draw  (v3) edge (v14);
\draw  (v4) edge (v15);
\draw  (v5) edge (v16);
\draw  (v5) edge (v17);
\draw  (v6) edge (v18);
\draw  (v6) edge (v19);

\node (v28) at (6,-6) {};
\node (v29) at (8,-6) {};
\node (v30) at (10,-6) {};
\node (v31) at (12,-6) {};
\draw  (8) edge (v28);
\draw  (8) edge (v29);
\draw  (v7) edge (v30);
\draw  (v7) edge (v31);

\end{tikzpicture}
\caption{Case 2.2.6.2.}
\label{2-case 2.2.6.2.}
\end{center}
\end{figure}
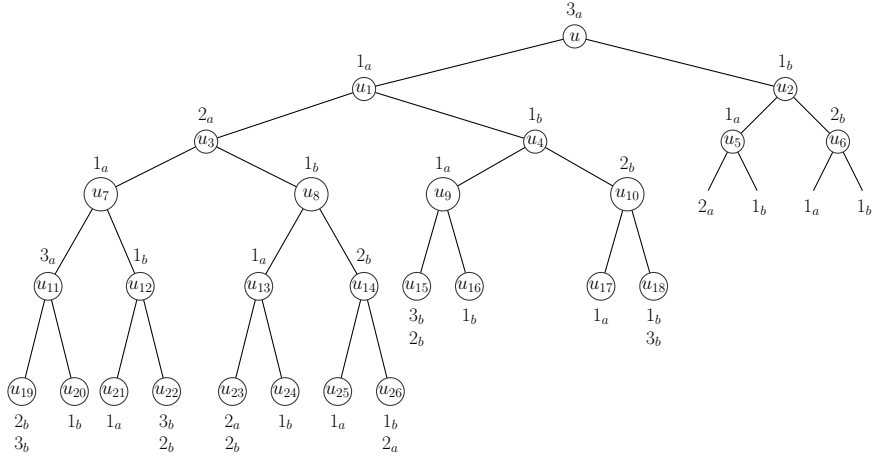


{\bf Case 2.2.6.2:} $f(u_{11})=3_a$ and $f(u_{14}) = 2_b$. Similarly to \eqref{u24u25}, we may assume
\begin{equation}\label{u20u21}
f(u_{20})=f(u_{24})=1_b \quad  \mbox{and}\quad   f(u_{21})=f(u_{25})=1_a. 
\end{equation}
Similarly to \eqref{2au13u14}, we may assume 
\begin{equation}\label{2au13u14'}
2_a \in f(N(u_{13}) \cup N(u_{14})-\{u_8\}). 
\end{equation}
If $1_b \notin f(N(u_{14})-\{u_8\})$ and $2_b \notin f(N(u_{13}))$, then we can recolor $u_8$ with $2_b$, $u_{14}$ with $1_b$, $u_3$ with $1_b$, and $u$ with $2_a$ to obtain a coloring satisfying $(a)\text{ and }(b)$. Thus, we may assume 
$$
f(N(u_{13})-\{u_8\})=\{1_b,2_a\}\quad  \mbox{and}\quad  f(N(u_{14})-\{u_8\})=\{1_a,1_b\}
$$
\begin{equation}\label{u13u14}
\quad \mbox{ or } \quad f(N(u_{13})-\{u_8\})=\{1_b,2_b\}\quad  \mbox{and}\quad  f(N(u_{14})-\{u_8\})=\{1_a,2_a\}.
\end{equation}
If $2_b \notin f(N(u_{11}) \cup N(u_{12}))$, then we can recolor $u_7$ with $2_b$ and it contradicts \eqref{cases}. If $3_b \notin f(N(u_{11}) \cup N(u_{12}))$, then we can recolor $u_3$ with $3_b$ and $u$ with $2_a$ to obtain a coloring satisfying $(a)\text{ and }(b)$. Thus, we may assume
\begin{equation}\label{u11u12'}
 f(N(u_{11}) \cup N(u_{12})-\{u_7\})=\{1_a,1_b,2_b,3_b\}.
\end{equation}
Specifically, we know that $1_a \notin f(N(u_{11})-\{u_7\})$ and $2_a \notin f(B(u_7,2) - \{u_3\})$ (See Figure~\ref{2-case 2.2.6.2.}). Therefore, we recolor $u_{11}$ with $1_a$, $u_7$ with $2_a$, $u_3$ with $3_a$, and $u$ with $2_a$ to obtain a coloring satisfying $(a)\text{ and }(b)$.  \hfill \qed

\bigskip\noindent
{\bf Acknowledgment.} We thank Sandi Klav\v zar, Douglas West, and the referees for their helpful comments.






\end{document}